%% file: GeneralCumulant3.tex
\documentclass[12pt,twoside]{mamsart}
\input{macros.tex}
\newcommand{\XX}{\mathfrak{X}} 
\newcommand{\bb}{\delta} 
\newcommand{\bx}{\boldsymbol{x}} 
\newcommand{\by}{\boldsymbol{y}} 
\newcommand{\dG}{\mathrm{d}\Gamma} 

\allowdisplaybreaks{}

\newcommand{\externalref}[2]{#1.\ref{#1-#2}}
\newcommand{\externaleqref}[2]{(#1.\ref{#1-#2})}
\newcommand{\ssection}[1]{\section{#1}}

\begin{document}

\title[Cumulants in Noncommutative Probability Theory III]
{Cumulants in Noncommutative Probability Theory III. Creation and
  annihilation operators on Fock spaces}
\author{Franz Lehner}

\thanks{Supported by the European Network \no{}HPRN-CT-2000-00116
and the Austrian Science Fund (FWF) Project \no{}R2-MAT}

\address{
Franz Lehner\\
In\-sti\-tut f\"ur Mathe\-ma\-tik C\\
Tech\-ni\-sche Uni\-ver\-si\-t\"at Graz\\
Stey\-rer\-gas\-se 30, A-8010 Graz\\
Austria}
\email{lehner@finanz.math.tu-graz.ac.at}


\date{\today}

\begin{abstract}
  Cumulants of noncommutative random variables arising from Fock space constructions
  are considered. In particular, simplified calculations are given for several known examples
  on~$q$-Fock spaces.

  In the second half of the paper we consider in detail the Fock states
  associated to characters of the infinite symmetric group recently constructed by
  Bo\.zejko and Guta. We express moments of multidimensional Dyck words in terms of
  the so called cycle indicator polynomials of certain digraphs.
\end{abstract}

\maketitle{}
\ddate{27.11.2000}

\tableofcontents{}

\input{Intro.tex}
\input{Prelim.tex}
\input{GeneralCumulantFock.tex}

\input{qFock.tex}
\newpage{}
\input{SymmetricGroup.tex}

\bibliography{GeneralCumulants}
\bibliographystyle{mamsalpha}
\end{document}

%% file: macros.tex
\usepackage[latin1]{inputenc}
\usepackage{amsmath,amsfonts,amssymb,amscd,amsthm,amsxtra}

\usepackage{a4wide}
\usepackage{url}
\usepackage{enumerate,verbatim}
\usepackage{ifthen}

\usepackage{url,textcomp}
\newcommand{\no}{\textnumero}

\usepackage{xr1}

\usepackage{pstricks}
\newcommand{\ddate}[1]{}

\newcommand{\ppaper}[4][ ]{}
\newcommand{\paper}[4][ ]{}

\numberwithin{equation}{section}

\makeatletter
  \DeclareRobustCommand\em
         {\@nomath\em \ifdim \fontdimen\@ne\font >\z@
                       \upshape \else \slshape \fi}
\makeatother

\newtheoremstyle{mythm}
  {9pt}
  {9pt}
  {\slshape}
  {0pt}
  {\bfseries}
  {.}
  { }
  {\thmname{#1} \thmnumber{#2}\thmnote{ (#3)}}

\theoremstyle{mythm} 
\newtheorem{Theorem}{Theorem}[section]
\newtheorem{Lemma}[Theorem]{Lemma}
\newtheorem{Claim}[Theorem]{Claim}
\newtheorem{Proposition}[Theorem]{Proposition}

\newtheorem{Corollary}[Theorem]{Corollary}

\theoremstyle{definition} 
\newtheorem{Definition}[Theorem]{Definition}

\newtheorem{Remark}[Theorem]{Remark}

\newtheorem{Example}[Theorem]{Example}
 
\newcommand{\file}[1]{{\url{#1}}}

\newcommand{\IC}{\mathbf{C}}                     
\newcommand{\IN}{\mathbf{N}}                     
\newcommand{\IR}{\mathbf{R}}                     
\newcommand{\IZ}{\mathbf{Z}}                     
\newcommand{\SG}{\mathfrak{S}}                   
\newcommand{\eps}{{\varepsilon}}                 
\newcommand{\oldphi}{{\phi}}                   
\renewcommand{\phi}{{\varphi}}                   
\newcommand{\ox}{\otimes}                        
\newcommand{\oxf}{\cdot}                        
\newcommand{\oxfdots}{\dotsm}                        
\newcommand{\NCPartitions}{NC}                     
\newcommand{\NC}{\NCPartitions}                    
\newcommand{\Fock}{\mathcal{F}}                    
\newcommand{\LCO}{L}                    
\newcommand{\RCO}{R}                    
\newcommand{\abs}[1]{\left\lvert #1 \right\rvert}  

\newcommand{\norm}[1]{\left\lVert #1 \right\rVert} 
\newcommand{\exchm}{\alg{E}} 
\newcommand{\exchF}{{\alg{F}}} 
\newcommand{\exch}{\mbox{$\exchm$}} 

\newlength{\tmpl}

\newcommand{\bra}[1]{\langle #1|}
\newcommand{\ket}[1]{|#1\rangle}

\newcommand{\indep}{\perp\kern-3pt\perp}

\newcommand{\alg}[1]{{\mathcal{#1}}}                
\newcommand{\hilb}[1]{{\mathfrak{#1}}}                

\DeclareMathOperator{\lin}{{\mathrm span}}                     



%% file: Intro.tex
\section*{Introduction}
Fock space constructions like in~\cite{GutaMaassen:2002:generalised}
give rise to natural interchangeable families
and are thus well suited for cumulant calculations like 
in part~I~\cite{Lehner:2002:Cumulants1} and
part~II~\cite{Lehner:2002:Cumulants2}.

In this paper we develop some general formulas and recompute cumulants
for generalized Toeplitz operators, notably for $q$-Fock spaces,
previously considered by A.~Nica~\cite{Nica_1996:crossings} 
and M.~Anshelevich~\cite{Anshelevich:2001:partitiondependent}, 
and the Fock spaces associated to characters of the infinite symmetric
group recently introduced by M.~Bo\.zejko and M.~Guta~\cite{BozejkoGuta:2002:functors}.
For the latter we indicate a general formula for mixed moments of creation and
annihilation operators in terms of the cycle indicator polynomial of a certain
directed graph.

The paper has three sections.

In the first section some general formulae are developed, in particular
the cumulants of generalized Toeplitz operators are given by
$$
K_n(L^* + \sum_{k=0}^\infty \alpha_{k+1} L^k)
= b_n\,\alpha_n
$$
where~$b_n$ is a certain statistic on the symmetric group.

In the second section, we review the $q$-cumulants of Nica and Anshelevich in
the light of the new theory.

In the final section we study in detail the Fock space associated to irreducible characters
on the infinite symmetric group recently introduced by Bo\.zejko and Guta.
It turns out that expectations of multidimensional Dyck words are given
by the so-called cycle-indicator polynomial of a certain digraph associated
to the Dyck word. This polynomial counts the number of coverings of the digraph
with hamiltonian cycles and satisfies a certain cut-and-fuse recursion formula,
which is reflected by a certain commutation relation of the creation and annihilation
operators. We characterize those Fock states, which only depend on this digraph construction,
as averages of the Bo\.zejko-Guta Fock states, that is, states which are
associated to not necessarily irreducible characters of the symmetric group.


%% file: Prelim.tex
\section{Preliminaries}
\label{sec:GeneralCumulant3:Preliminaries}

In this section we collect the necessary definitions and auxiliary results
needed later on. For details we refer to part I~\cite{Lehner:2002:Cumulants1}.

\subsection{Exchangeability Systems and Cumulants}

We recall first that an \emph{exchangeability system} $\exch$ for a
noncommutative probability space~$(\alg{A},\phi)$ consists of another
noncommutative probability space $(\alg{U},\tilde{\phi})$ and an infinite
family $\alg{J}=(\iota_k)_{k\in\IN}$ of state-preserving embeddings
$\iota_k:\alg{A}\to \alg{A}_k\subseteq \alg{U}$, 
which we conveniently denote by $X\mapsto X^{(k)}$,
such that the algebras $\alg{A}_j$ are \emph{interchangeable}
with respect to $\tilde\phi$:
for any family~$X_1,X_2,\dots,X_n\in \alg{A}$,
and for any choice of indices~$h(1),\dots,h(n)$
the expectation is invariant under any permutation~$\sigma\in\SG_\infty$
in the sense that
\begin{equation}
  \label{eq:GeneralCumulant4:phih=phisigmah}
  \tilde\phi(X_1^{(h(1))} X_2^{(h(2))} \cdots X_n^{(h(n))})
  = \tilde\phi(X_1^{(\sigma(h(1)))} X_2^{(\sigma(h(2)))} \cdots X_n^{(\sigma(h(n)))})
  .
\end{equation}
Throughout this paper we will assume that the algebra $\alg{U}$ is generated
by the algebras $\alg{A}_k$ and that the action of $\SG_\infty$ extends to
all of $\alg{U}$.
Denote by~$\Pi_n$ the set of partitions (or equivalence relations) of the
set $[n]=\{1,2,\dots,n\}$.
The value~\eqref{eq:GeneralCumulant4:phih=phisigmah} only depends on the
kernel~$\pi=\ker h\in\Pi_n$ defined by 
$$
i\sim_\pi j
\qquad
\iff
\qquad
h(i)=h(j)
$$
and we denote it
$
\phi_\pi(X_1,X_2,\dots,X_n) = 
\tilde\phi(X_1^{(\pi(1))} X_2^{(\pi(2))} \cdots X_n^{(\pi(n))})
.
$
Here we consider a partition~$\pi\in\Pi_n$ as a function $\pi:[n]\to\IN$,
mapping each element to the number of the block containing it.
This is a canonical example of an index function~$h$ with $\ker h=\pi$.

Subalgebras $\alg{B},\alg{C}\subseteq \alg{A}$ are called
\emph{$\alg{E}$-exchangeable} or, more suggestively, \emph{$\alg{E}$-independent}
if for any choice of random variables~$X_1,X_2,\dots,X_n \in \alg{B}\cup\alg{C}$
and subsets~$I,J\subseteq \{1,\dots,n\}$ such that~$I\cap J=\emptyset$,
$I\cup J=\{1,\dots,n\}$, $X_i\in \alg{B}$ for $i\in I$  and $X_i\in \alg{C}$ for $i\in J$,
we have the identity
$$
\phi_\pi(X_1,X_2,\dots,X_n)
= \phi_{\pi'}(X_1,X_2,\dots,X_n)
$$
whenever~$\pi$, $\pi'\in\Pi_n$ are partitions with~$\pi|_I=\pi'|_I$ and
$\pi|_J=\pi'|_J$.
We say that two families of random variables~$(X_i)$ and~$(Y_j)$ are
$\alg{E}$-exchangeable if the algebras they generate have this
property.


Then it is possible to define cumulant functionals, indexed by set partitions
$\pi\in\Pi_n$, via
$$
K_\pi^\exchm(X_1,X_2,\dots,X_n)
= \sum_{\sigma\leq\pi}
   \phi_\sigma(X_1,X_2,\dots,X_n)\,\mu(\sigma,\pi) 
$$
where $\mu(\sigma,\pi)$ is the M\"obius function of the lattice of set
partitions, cf.\ part~I.
Alternatively, the cumulants can be defined by Good's formula.
Given noncommutative random variables~$X_1$,~$X_2$, \ldots,~$X_n$,
take an arbitrary partition~$\pi\in\Pi_n$.
We choose for each $k\in \{1,\dots,n\}$ an exchangeable
copy~$\{X_j^{(k)}:j\in \{1,\dots,n\}\}$ of the given family~$\{X_j:j\in \{1,\dots,n\}\}$.
For each block~$B=\{k_1<k_2<\dots<k_b\}\in\pi$ we pick
a primitive root of unity~$\omega_b$ of order~$b=|B|$, 
and set for each~$i\in B$
$$
X_i^{\pi,\omega}
= \omega_b X_i^{(k_1)}
  +
  \omega_b^2 X_i^{(k_2)}
  +
  \dots
  +
  \omega_b^{b} X_i^{(k_b)}
$$
Then we have Good's formula~\cite[Prop.~2.8]{Lehner:2002:Cumulants1}
\begin{equation}
  \label{eq:GeneralCumulants:PartitionGoodFormula}
  K^\exchm_\pi(X_1,X_2,\dots,X_n)= \frac{1}%
  {\prod \abs{B}}
  \,
  \phi(X_1^{\pi,\omega}X_2^{\pi,\omega}\dotsm X_n^{\pi,\omega})
  ;
\end{equation}
if~$\pi\in\Pi_n$ consists of one block only, we may abbreviate and write
$$
K_n(X_1,X_2,\dots,X_n) = \frac{1}{n}\, \phi(X_1^{\omega}X_2^{\omega}\dotsm X_n^{\omega})
$$
where~$\omega$ is a primitive root of unity of order~$n$ and
$$
X_i^\omega
= \omega X_i^{(1)}
  +
  \omega^2 X_i^{(2)}
  +
  \dots
  +
  \omega^n X_i^{(n)}
  .
$$

The use of the probabilistic termini ``independence'' and ``cumulants'' is
justified by the following proposition which establishes the analogy to
classical probability.

\begin{Proposition}[{\cite{Lehner:2002:Cumulants1}}]
  Two subalgebras~$\alg{B},\alg{C}\subseteq\alg{A}$ are \exch-independent 
  if and only if mixed cumulants vanish, that is,
  whenever~$X_i\in\alg{B}\cup\alg{C}$ are some noncommutative random variables
  and~$\pi\in\Pi_n$ is an arbitrary partition
  such that there is a block of $\pi$ which contains indices $i$ and $j$
  such that~$X_i\in\alg{B}$ and~$X_j\in\alg{C}$,
  then $K^\exchm_\pi(X_1,X_2,\dots,X_n)$ vanishes.
\end{Proposition}

\begin{Remark}
\label{rem:GeneralCumulant4:iid=interchangeable}
  In the sequel we will sometimes not distinguish between i.i.d.\ sequences
  in~$\alg{A}$ and sequences of the form~$X^{(i)}$. The latter do not belong to~$\alg{A}$
  strictly speaking, but we can replace $\alg{A}$ by the algebra~$\tilde{\alg{A}}$
  generated by~$(\alg{A}_i)_{i\in I}$, where~$I\subseteq\IN$ is an infinite
  subset, and construct an exchangeability system for~$\tilde{\alg{A}}$ by
  considering~$\IN$ as a disjoint union of infinitely many copies of~$I$.
\end{Remark}
One of our main tools will be the product formula of Leonov and Shiryaev.
\begin{Proposition}[{\cite[Prop.~3.3]{Lehner:2002:Cumulants1}}]
  \label{prop:GeneralCumulants:Productformula}
  Let~$(X_{i,j})_{i\in\{1,\dots,m\},j\in \{1,\dots,n_i\}}\subseteq\alg{A}$
  be a family of noncommutative random variables,
  in total $n=n_1+n_2+\dots+n_m$ variables.
  Then every partition~$\pi\in\Pi_m$ induces a partition~$\tilde{\pi}$
  on $\{1,\dots,n\}\simeq\{(i,j):i\in [m], j\in[n_i]\}$ with blocks
  $\tilde{B}=\{(i,j):i\in B, j\in[n_i]\}$, that is, each block~$B\in\pi$ is
  replaced by the union of the intervals~$(\{n_{i-1}+1,n_{i-1}+2,\dots,n_i\})_{i\in B}$.
  Then we have
  $$
  K^\exchm_\pi(\prod_{j_1} X_{1,j_1},\prod_{j_2} X_{2,j_2},\dots,\prod_{j_m} X_{m,j_m})
  =\sum_{\substack{\sigma\in\Pi_n\\ \sigma\vee \tilde{\hat0}_m=\tilde\pi}}
    K^\exchm_\sigma(X_{1,1},X_{1,2},\dots,X_{m,n_m})
  $$
\end{Proposition}


%% file: GeneralCumulantFock.tex
\section{Fock spaces}
\ddate{15.02.2002}
We recall some definitions from~\cite{GutaMaassen:2002:generalised}.
\begin{Definition}
  Let $\hilb{K}$ be a real Hilbert space.
  The algebra $\hilb{A}(\hilb{K})$ is the unital $*$-algebra with
  generators $\{\omega(\xi):\xi\in\hilb{K}\}$ and relations
  $$
  \omega( \lambda \xi+\mu \eta) = \lambda\,\omega(\xi)+\mu\,\omega(\eta)
  \qquad
  \omega(\xi)^*=\omega(\xi)
  $$
  for all~$\xi$,~$\eta\in\hilb{K}$ and~$\lambda$,~$\mu\in\IR$.
\end{Definition}
\begin{Definition}
  Let $\hilb{H}$ be a complex Hilbert space.
  The algebra $\alg{C}(\hilb{H})$ is the unital $*$-algebra with generators
  $\LCO(\xi)$ and $\LCO^*(\xi)$ and relations
  $$
  \LCO(\lambda \xi+\mu \eta) = \lambda \LCO(\xi) +\mu \LCO(\eta)
  \qquad
  \LCO^*(\xi) = \LCO(\xi)^*
  $$
  for all~$\xi$,~$\eta\in\hilb{H}$ and~$\lambda$,~$\mu\in\IC$.
  
\end{Definition}

\begin{Definition}
  A \emph{Fock state} on $\alg{C}(\hilb{H})$ is a state $\rho$ satisfying
  $$
  \begin{bmatrix}
    \rho(\LCO^*(\xi)\,\LCO^*(\eta)) &     \rho(\LCO^*(\xi)\,\LCO(\eta)) \\
    \rho(\LCO(\xi)\,\LCO^*(\eta)) &       \rho(\LCO(\xi)\,\LCO(\eta))
  \end{bmatrix}
  =
  \begin{bmatrix}
    0 & \langle \xi,\eta\rangle \\
    0 & 0
  \end{bmatrix}
  $$
  and 
  $$
  \rho(\LCO^{\eps_1}(U\xi_1)\, \LCO^{\eps_2}(U\xi_2) \cdots \LCO^{\eps_n}(U\xi_n))
  = \rho(\LCO^{\eps_1}(\xi_1)\, \LCO^{\eps_2}(\xi_2)\cdots \LCO^{\eps_n}(\xi_n))
  $$
  for all unitary operators $U\in\alg{U}(\hilb{H})$
  and all choices of exponents $\eps_j\in\{*,1\}$.
  In particular, odd moments vanish.
\end{Definition}

If the Hilbert space $\hilb{H}$ is infinite dimensional,
we can decompose it into an infinite direct sum
$\hilb{H}=\bigoplus \hilb{H}_i$ and the subalgebras
$\alg{C}(\hilb{H}_i)$ are interchangeable, because any permutation
of indices can be implemented by a unitary operator on $\hilb{H}$.
Therefore Fock spaces are a rich source of exchangeability systems.
We proceed by calculating cumulants in certain special cases.
Throughout this paper the exchangeability system will be
$$
\exch=(\alg{C}(\hilb{K}),\rho, (\iota_k))
$$
where $\hilb{K}=\bigoplus\hilb{H}_i$ is the direct sum of infinitely many
copies $\hilb{H}_i\simeq\hilb{H}$, which give rise
to the embeddings $\iota_k : \alg{C}(\hilb{H})\to\alg{C}(\hilb{H}_k)\subseteq\alg{C}(\hilb{K})$,
which are interchangeable with respect to the Fock state $\rho$.

\begin{Proposition}
  \label{prop:GeneralCumulants:FockState}
  Any Fock state is given by a function $\mathbf{t}$ on pair partitions
  \begin{equation}
    \label{eq:GeneralCumulants:FockState=sumPi2}
  \rho(\LCO^{\eps_1}(\xi_1)\,\LCO^{\eps_2}(\xi_2)\cdots \LCO^{\eps_n}(\xi_n))
  = \sum_{\pi\in\Pi_n^{(2)}} \mathbf{t}(\pi)
     \prod_{\{k<l\}\in\pi}
     Q(\eps_k,\eps_l)
     \,
     \langle \xi_k, \xi_l\rangle
  \end{equation}
  where $\eps_j\in\{*,1\}$ and
  $$
  Q(\eps,\eps') = 
  \begin{cases}
    1 & \text{if $\eps=*$ and $\eps'=1$}\\
    0 & \text{otherwise}
  \end{cases}
  $$
\end{Proposition}
The proof is essentially the same as the proof of Theorem~\externalref{II}{thm:GeneralCumulants:Maxwell}.
Conversely, it was shown in~\cite{GutaMaassen:2002:generalised}
that a positive definite function on pair partitions gives rise to a Fock state
and a Fock representation of $\alg{C}(\hilb{H})$.

\begin{Corollary}
  For $\eps_j\in\{1,*\}$
  the cumulants of the creation and annihilation operators are given by
  $$
  K^\exchm_\pi(\LCO^{\eps_1}(\xi_1), 
        \LCO^{\eps_2}(\xi_2),
        \cdots,
        \LCO^{\eps_n}(\xi_n))
  = \begin{cases}
      0 & \text{if $\pi\not\in\Pi_n^{(2)}$} \\
     \mathbf{t}(\pi)
     \prod_{\{k<l\}\in\pi}
     Q(\eps_k,\eps_l)
     \,
     \langle \xi_k, \xi_l\rangle
        & \text{if $\pi\in\Pi_n^{(2)}$}
    \end{cases}
  $$
\end{Corollary}

\begin{Definition}
  A positive definite function $\mathbf{t}$ on pair partitions is called
  \emph{multiplicative} if it factors with respect to the connected components of $\pi$.
  (See Definition~\externalref{I}{def:GeneralCumulants:partitions}).
\end{Definition}
The following lemma is immediate from the definition.
\begin{Lemma}
  Pyramidal independence holds if and only if~$\mathbf{t}$ is multiplicative.
\end{Lemma}
From now on we will work with a fixed orthonormal basis $\{e_i\}$ of $\hilb{H}$
and calculate expectation of words in $\LCO_i=\LCO(e_i)$ or creation
operators $\LCO=\LCO(h)$ for some fixed unit vector $h\in\hilb{H}$.
\begin{Definition}
  \label{def:GeneralCumulant3:latticepath}
  A \emph{lattice path} is a sequence of points
  $((x_i,y_i))_{i=0,1,\dots,n}$ in $\IN_0\times\IZ$
  such that $y_0=y_n=0$, $y_j\geq0$ for all $j$ and $x_i=i$.
  As the $x$-coordinates are redundant, 
  we will also refer to the sequence $(y_i)_{i=0,1,\dots,n}$ as lattice paths.
  A lattice path is \emph{irreducible} if $y_j>0$ for $j=1,2,\dots,n-1$
  and reducible otherwise.
  A \emph{\L{}ukasiewicz path} is a lattice path $(y_j)$ such that
  $y_i-y_{i-1}\leq1$.
  A \emph{lattice word} is a word $\LCO^{k_1}\LCO^{k_2}\cdots \LCO^{k_n}$
  where $k_j\in\IZ$ and $\LCO^{-k}$ is interpreted as $\LCO^*{}^k$ and 
  such that the sequence $y_j = k_1+k_2+\dots+k_j$ constitutes a lattice path.
  A \L{}ukasiewicz word is defined accordingly.
  A \emph{Dyck path} is a lattice path $(y_j)$ such that $y_i-y_{i-1}=\pm 1$.
  A \emph{$N$-dimensional lattice path} is a sequence of points 
  $\vec{y}_j$ in $\IZ^N$ such that $\vec{y}_0=\vec{y}_n=\vec{0}$
  and $\vec{y}_{k+1}-\vec{y}_{k}$ is a multiple of a basis vector
  $e_i$ and such that all coordinates of $\vec{y}_j$ are nonnegative.
  A word $\LCO_{i_1}^{k_1}\LCO_{i_2}^{k_2}\cdots \LCO_{i_n}^{k_n}$ is an
  \emph{($N$-dimensional) lattice word}
  if the path $\vec{y}_{j}=k_1 e_{i_1} + k_2 e_{i_2} + \dots + k_j e_{i_j}$
  is a lattice path. Multidimensional \L{}ukasiewicz words and Dyck words are defined
  similarly.
\end{Definition}
\begin{Proposition}
  For $k_j\in\IZ$ we have
  $$
  \rho(\LCO_{i_1}^{k_1}\LCO_{i_2}^{k_2}\cdots \LCO_{i_n}^{k_n}) = 0
  $$
  unless the word $\LCO_{i_1}^{k_1}\LCO_{i_2}^{k_2}\cdots \LCO_{i_n}^{k_n}$
  is a lattice word. 
\end{Proposition}
\begin{proof}
  Indeed if one of the $\vec{y}_j$ in Definition~\ref{def:GeneralCumulant3:latticepath} 
  has a negative component
  then the sum \eqref{eq:GeneralCumulants:FockState=sumPi2} is empty,
  as there is no pairing with nonzero contribution. 
\end{proof}

\ddate{17.07.2002}
It was shown in~\cite{GutaMaassen:2002:generalised} that Fock states can
always be modeled on so-called combinatorial Fock spaces as follows.
Let $\hilb{H}$ be a fixed Hilbert space.
Let $V_n$ be a sequence of Hilbert spaces with an action $U_n$ of
the symmetric group $\SG_n$ and densely defined intertwining operators
$j_n:V_n\to V_{n+1}$ such that
\begin{equation}
  \label{eq:GeneralCumulantFock:intertwining}
  j_n\circ U_n(\sigma) = U_{n+1}(\iota_n(\sigma))\circ j_n
\end{equation}
where $\iota_n:\SG_n\to\SG_{n+1}$ is the natural inclusion.
Let
$$
\Fock_V(\hilb{H})
= \bigoplus_{n=0}^\infty \frac{1}{n!} \,V_n\ox_s\hilb{H}^{\ox n}
$$
where
$V_n\ox_s\hilb{H}^{\ox n}$ is the subspace spanned by the vectors which are invariant
under the action $U_n\ox\tilde{U}_n$ of $\SG_n$ on $V_n\ox \hilb{H}^{\ox n}$ given by
$$
(U_n(\sigma)\ox\tilde{U}_n(\sigma))\, (v\ox\xi_1\ox \xi_2\ox\cdots \ox \xi_n)
= U_n(\sigma)\,
  v\ox
  \xi_{\sigma^{-1}(1)}
  \ox
  \xi_{\sigma^{-1}(2)}
  \ox 
  \cdots
  \ox
  \xi_{\sigma^{-1}(n)}
$$
From now on we will abbreviate
$$
v\ox\xi_1\cdot \xi_2\oxfdots  \xi_n:=v\ox\xi_1\ox \xi_2\ox\cdots \ox \xi_n
;
$$
for~$v\in V_n$ and~$\eta\in\hilb{H}^{\ox n}$ we define the Symmetrizator
$$
v\ox_s \eta
= P_n v\ox\eta
= \frac{1}{n!}
  \sum_\sigma 
   U_n(\sigma)\, v\ox \tilde{U}_n\,(\sigma)\eta
$$
and left and right creation operators
\begin{equation}
  \label{eq:GeneralCumulantFock:creationoperatordefinition}
  \LCO_{V,j}(h) \, v\ox_s\eta
  = (n+1)\, (j_nv)\ox_s (h\ox\eta)
  \qquad\qquad
  \RCO_{V,j}(h)\, v\ox_s\eta
  = (n+1)\, (j_nv)\ox_s (\eta\ox h)
  ;
\end{equation}
its adjoint is the annihilation operator $\LCO^*_{V,j}$ which is the restriction of
\begin{align*}
  \tilde{\LCO}^*_{V,j}(\xi):V_{n+1}\ox\hilb{H}^{\ox n+1}& \to V_n\ox\hilb{H}^{\ox n}\\
   v\ox \xi_1\oxf \xi_2\oxfdots \xi_n &\mapsto \langle \xi_1,\xi\rangle\, 
                                            j_n^* v\ox \xi_2\cdot\xi_3\oxfdots \xi_n
\end{align*}
to $V_{n+1}\ox_s \hilb{H}^{\ox n+1}$.
A similar formula holds for $\RCO^*_{V,j}$. On symmetric tensors it is given by
\begin{equation}
  \label{eq:GeneralCumulantFock:annihilationoperator}
  \RCO^*(\xi) \, v\ox_s \xi_1\oxfdots\xi_{n+1}
  =\frac{1}{n+1} 
   \sum_{i=1}^{n+1} \langle \xi_i,\xi \rangle \,
    j_n^*U_{n+1}(\tau_{i,n+1})\, v \ox_s \xi_1\oxfdots\xi_{i-1}\oxf\xi_{n+1}\oxf\xi_{i+1}\oxfdots\xi_n
\end{equation}
where $\tau_{i,j}$ is the transposition which exchanges $i$ and $j$.
Left and right creation operators are equivalent and we will consider either of them,
whenever it is notationally convenient.

Let $\alg{C}_{V,j}(\hilb{H})$ be the $*$-algebra generated by these creation operators
and $\Omega_V\in V_0$ a unit vector.
Then
$$
\rho_{V,j}(X) = \langle \Omega_V,X\Omega_V \rangle
$$
is a Fock state.
\begin{Proposition}{{\cite[Thm.~2.7]{GutaMaassen:2002:generalised}}}
  Let $\mathbf{t}$ be a positive definite function on pair partitions.
  Then for any complex Hilbert space~$\hilb{H}$
  the GNS-representation of $\alg{C}(\hilb{H},\rho_{\mathbf{t}})$
  is unitarily equivalent to $(\Fock_V(\hilb{H}),\alg{C}_{V,j}(\hilb{H}),\Omega_V)$ for some
  sequence $(V_n,j_n)_{n=0}^\infty$ which up to unitary equivalence is uniquely
  determined by $\mathbf{t}$.
\end{Proposition}

Using this model it is now easy to compute cumulants in certain cases.
\begin{Proposition}
  If the lattice word $(\LCO^{k_1},\LCO^{k_2}, \dots,\LCO^{k_n})$ is reducible,
  then 
  $$
  K^\exchm_n(\LCO^{k_1},\LCO^{k_2},\dots,\LCO^{k_n})=0
  $$
\end{Proposition}
\begin{proof}
  Indeed 
  for a lattice word $\LCO_{i_1}^{k_1}\LCO_{i_2}^{k_2}\cdots \LCO_{i_m}^{k_m}$
  we have
  $$
  \LCO_{i_1}^{k_1}\LCO_{i_2}^{k_2}\cdots \LCO_{i_m}^{k_m} \Omega_V
  = \rho(\LCO_{i_1}^{k_1}\LCO_{i_2}^{k_2}\cdots \LCO_{i_m}^{k_m})
    \,
    \Omega_V
  $$
  and therefore
  if $m$ is the largest index for which $k_{m+1}+k_{m+2}+\dots+k_n=0$, then
  by assumption $m>1$ and
  Lemma~\externalref{I}{lem:GeneralCumulants:subwordexpectationvanishes}
  implies that
  $$
  (\LCO^{k_{m+1}})^\omega (\LCO^{k_2})^\omega \dots (\LCO^{k_n})^\omega \Omega_V
  = \rho((\LCO^{k_{m+1}})^\omega (\LCO^{k_2})^\omega \dots (\LCO^{k_n})^\omega)
    \,
    \Omega_V
  = 0
  .
  $$
\end{proof}

In the case of the free creation operators of Voiculescu
we can say even more, 
see Proposition~\ref{prop:GeneralCumulants:FreeFockCumulants} below.
There is a special kind of operators for which cumulants can be calculated
explicitly in general.
\begin{Definition}
  Let $h\in\hilb{H}$ be a fixed unit vector and let $\LCO=\LCO(h)$ be a creation
  operator. 
  A \emph{generalized Toeplitz operator} is an operator of the form
  $$
  \LCO^* + \sum_{k=1}^\infty \alpha_k \LCO^{k-1}
  $$
\end{Definition}
Such operators were considered by Voiculescu as a model of free random variables
\cite{Voiculescu:1986:addition}
and for $q$-deformations of them by Nica~\cite{Nica_1996:crossings}.
The cumulants of Toeplitz operators are rather restrictive and only allow
\L{}ukasiewicz words which are ``prime'' and consequently can be immediately read off
the coefficients.
\begin{Proposition}
  \label{prop:GeneralCumulantFock:Lukasiewiczcumulants}
  For $k_j\in\{*,0,1,\dots\}$ we have
  $$
  K^\exchm_n(\LCO^{k_1},\LCO^{k_2},\dots,\LCO^{k_n}) = 0
  $$
  unless $k_1=k_2=\cdots=k_{n-1}=*$ and $k_n=n-1$.
  In the latter case, the cumulant equals the expectation of the product:
  \begin{align*}
    K^\exchm_n(\LCO^*,\dots,\LCO^*,\LCO^{n-1})
    &= \rho((\LCO^*)^{n-1} \LCO^{n-1}) \\
    &= \sum_{\pi\in\Pi_{2n-2}^{(2)}}
        K^\exchm_\pi(\LCO^*,\LCO^*,\dots,\LCO,\LCO)\\
    &= \sum_{\sigma\in\SG_{n-1}} \nu(\sigma)
  \end{align*}
  where $\nu$ is a function on the symmetric group,
  because every contributing $2$-partition connects a point of
  $\{1,2,\dots,n-1\}$ to a point in $\{n,n+1,\dots,2n-2\}$
  and can be interpreted as a permutation.
\end{Proposition}
\begin{proof}
  We may assume that all $k_j\ne0$,
  because otherwise the presence of the identity operator makes the cumulant vanish.
  Let $m$ be the total number of factors when decomposing the powers of $\LCO$
  into single creators, i.e. $m=\sum\abs{k_j}$ where we put $\abs{*}=1$.
  Let $\pi\in\Pi_m$ be the interval partition induced by the $n$ exponents $k_j$,
  that is $\pi=\tilde{\hat0}_n$ in the notation of the proof
  of Proposition~\externalref{I}{prop:GeneralCumulants:Productformula}.
  If the $k_j$ are not as claimed, then there are at least two monomials $\LCO^k$
  with $k>0$.
  By the product formula (Proposition~\externalref{I}{prop:GeneralCumulants:Productformula})
  we have
  $$
  K^\exchm_n(\LCO^{k_1},\LCO^{k_2},\dots,\LCO^{k_n})
  = \sum_{\sigma\vee\pi=\hat1_m}
    K^\exchm_\sigma(\LCO^{\eps_1},\LCO^{\eps_2},\dots,\LCO^{\eps_m})
  $$
  where $\eps_j\in\{*,1\}$. Now $\sigma$ must be both a pair partition 
  connecting creation operators with annihilation operators and at the 
  same time we must have $\sigma\vee\pi=\hat1_m$,
  but the two different blocks with~$k>0$
  cannot be connected in this way, because each $L^*$ can only be connected to
  one $L^k$. Therefore the sum is empty.
\end{proof}

A similar formula holds for partitioned cumulants.
\begin{Proposition}
  \label{prop:GeneralCumulantFock:PartitionedLukasiewiczcumulants}
  For $k_j\in\{*,0,1,\dots\}$ we have
  $$
  K^\exchm_\pi(\LCO^{k_1},\LCO^{k_2},\dots,\LCO^{k_n}) = 0
  $$
  unless $k_j$ and $\pi$ are compatible in the sense that
  $$
  k_j
  =
  \begin{cases}
    b-1 & \text{if $j$ is the last element of a block $B$ of length
    $\abs{B}=b$}\\
    *   & \text{otherwise}
  \end{cases}
  $$
  In that case 
  $$
  K^\exchm_\pi(\LCO^{k_1},\LCO^{k_2},\dots,\LCO^{k_n})
  = \rho_\pi(\LCO^{k_1},\LCO^{k_2},\dots,\LCO^{k_n})
  $$
\end{Proposition}

By multilinear expansion we have the following corollary.
\begin{Corollary}
  \label{cor:GeneralCumulantFock:generalizedToeplitzcumulant}
  Let
  $$
  b_n
  = \rho((\LCO^*)^{n-1} \LCO^{n-1}) 
  = \sum_{\sigma\in\SG_{n-1}} \nu(\sigma)
  $$
  where $\nu(\sigma)$ is the statistic on the symmetric group defined in 
  Proposition~\ref{prop:GeneralCumulantFock:Lukasiewiczcumulants}.
  Then the generalized Toeplitz operator
  $$
  T = \LCO^* + \sum_{k=1}^\infty \frac{\alpha_k}{b_k} \LCO^{k-1}
  $$
  has cumulants
  $$
  K^\exchm_n(T,T,\dots,T) = \alpha_n
  $$
\end{Corollary}

\begin{Example}[$q$-cumulants]
  For Nica's $q$-cumulants \cite{Nica_1996:crossings} we have 
  $$
  b_n(q) = [n]_q!
  $$
  which includes the classical case $b_n(1)=n!$
  and the free case $b_n(0)=1$ considered by Voiculescu.
  Moreover, for partitioned cumulants a simple formula holds, too,
  namely Corollary~\ref{cor:qCumulantsN:partitionedcumulants} below.
\end{Example}


%% file: qFock.tex
\goodbreak{}
\section{$q$-Fock space}

\ddate{20.07.2001}
In an attempt to unify bosonic and fermionic Fock space,
the deformed $q$-Fock spaces were constructed in 
\cite{FrischBourret:1970:parastochastics}%
\ppaper{Frisch/Bourret}{Parastochastics}{J. Mathematical Phys. 11 (1970)
  364--390},
\cite{BozejkoSpeicher:1994:completely}%
\ppaper{Bozejko/ Speicher}{Completely positive maps on {C}oxeter groups, deformed commutation relations, and operator spaces}{Math. Ann. 300 (1994) 97--120}
and
\cite{BozejkoKuemmererSpeicher:1997:qGaussian}%
\ppaper{Bozejko/Kümmerer/Speicher}{$q$-Gaussian processes: non-commutative and
  classical aspects}{Comm. Math. Phys. 185 (1997) 129--154}.

\begin{Definition}[$q$-Fock space]
  On free Fock space 
  $\Fock(\hilb{H}) = \Fock_0(\hilb{H}) = \IC\Omega\oplus\bigoplus_{n\geq 1} \hilb{H}^{\ox n}$
  define the $q$-symmetrizator on $n$-particle space by
  \begin{align*}
    P_n^q(\eta_1\oxf \eta_2\oxfdots\eta_n)
    &= \sum_{\sigma\in\SG_n}
        q^{\abs{\sigma}} 
        \tilde{U}_n(\sigma)\,\eta_1\oxf \eta_2\oxfdots\eta_n \\
    &= \sum_{\sigma\in\SG_n}
        q^{\abs{\sigma}} \,
        \eta_{\sigma^{-1}(1)}\oxf \eta_{\sigma^{-1}(2)}\oxfdots\eta_{\sigma^{-1}(n)}
  \end{align*}
  where $\abs{\sigma}$ is the number of inversions of the partition $\sigma$.
  Define the $q$-inner product on elementary tensors $\xi\in \hilb{H}^{\ox m}$,
  $\eta\in \hilb{H}^{\ox n}$ by
  $$
  \langle \xi,\eta\rangle_q = \delta_{m,n}\, \langle \xi, P_n^q \eta\rangle
  $$
  where $\langle\,,\,\rangle$ is the inner product on $\Fock(\hilb{H})$ (linear in the second variable).
  $q$-Fock space is the completion with respect to this norm and denoted $\Fock_q(\hilb{H})$.
  On this space we define creation operators and annihilation operators
  \begin{align}
    \LCO(\xi)\,\Omega &= \xi 
    & \LCO(\xi)\,\eta_1\oxfdots\eta_n &=\xi\oxf\eta_1\oxfdots\eta_n \\
    \label{eq:GeneralCumulantExamples:qFock:annihilation}
    \LCO^*(\xi)\, \Omega &= 0 
    & \LCO^*(\xi)\,\eta_1\oxfdots\eta_n
        &=\sum q^{k-1}
               \langle\xi,\eta_k\rangle
               \,
               \eta_1 \oxfdots \hat{\eta}_k \oxfdots \eta_n
  \end{align}
\end{Definition}
  These satisfy the $q$-commutation relations
  $$
  \LCO^*(\xi)\,\LCO(\eta) - q\, \LCO(\eta)\,\LCO^*(\xi) = \langle \xi,\eta\rangle\, I
  $$
  and give rise to canonical interchangeable random variables
  with respect to the vacuum expectation
  $\rho(X) = \langle \Omega, X\Omega\rangle$:
  Given mutually orthogonal and isomorphic subspaces $(E_j)_{j\geq0}$
  of $\hilb{H}$, the $*$-subalgebras $\alg{A}_j$ generated by $\{\LCO(f):f\in E_j\}$
  are interchangeable. Moreover, orthogonal subspaces of $E_0$ give
  rise to interchangeable algebras.
  \begin{Proposition}[{\cite{BozejkoSpeicher:1994:completely}}]
    The Fock state corresponding to $q$-Fock space is given by the 
    positive definite function
    $$
    \mathbf{t}(\pi) = q^{nc(\pi)}
    $$
    where $nc(\pi)$ is the \emph{number of crossings} of the pair
    partition~$\pi$.
    In particular, pyramidal independence holds.
  \end{Proposition}
  The cases $q=-1,0,1$ give rise to fermionic, free and bosonic
  Fock spaces and correspondingly fermionic graded probability theory
  \cite{MingoNica:1997:crossings}%
  \paper{Mingo/Nica}{}{} (see section~\externalref{I}{sec:GeneralCumulantExamples:Fermions}),
  free probability theory \cite{Voiculescu:1985:symmetries,Voiculescu:1986:addition}
  and classical probability theory with corresponding notions of independence and convolution.
  
  However, for other $q$ there is no $q$-convolution
  \cite{vanLeeuwenMaassen:1996:obstruction}%
  \paper{van Leeuwen/Maassen}{}{}:
  There exist $q$-independent s.a.\ variables $X$, $X'$, $Y$ such that
  $X$ and $X'$ have the same distribution but the distributions
  of $X+Y$ and $X'+Y$ differ.
  See also \cite{Guta:2001:qproduct}, where a $q$-convolution for generalized
  Gaussians is constructed.

  Nevertheless, choosing certain models of random variables
  various $q$-cumulants have been found.

\ddate{25.07.2001}
\subsection{$q$-Toeplitz operators}
Motivated by Speicher's work on free cumulants and Voiculescu's
Toeplitz model of free random variables~\cite{Voiculescu:1986:addition}, Nica 
\cite{Nica:1995:oneparameter}%
\paper{Nica}{A one-parameter family of transforms, linearizing convolution laws for probability distributions}{Comm. Math. Phys. 168 (1995) 187--207}
considered $q$-Toeplitz operators 
\begin{equation}
  \label{eq:qCumulantsN:ToeplitzModel}
  T=\LCO_q^* + \sum_{n=0}^\infty \frac{ \alpha_{k+1}}{[k]_q!}\, \LCO_q^k
\end{equation}
with $L=L(h)$ for some fixed unit vector $h$.
Here and in the sequel we will use the standard $q$-notations
$$
  [n]_q = \frac{1-q^n}{1-q} = 1+q+q^2+\dots+q^{n-1}
  \qquad  \qquad
  [n]_q! = [1]_q [2]_q \cdots [n]_q
$$
Here exchangeable copies correspond to
$\LCO_i = \LCO_q^{(i)} = \LCO_q(e_i)$.
We want to show that Nica's cumulants
involving left-reduced crossings coincide with the cumulants of this exchangeability system.
Actually technical reasons force us to rather consider right-reduced crossings.
We will calculate mixed cumulants of creation- and annihilation operators.




\begin{Theorem}
  For a lattice word $(\eps_1,\eps_2,\dots,\eps_n)$,
  $\eps_j\in\{*,0,1,\dots\}$,
  the $q$-cumulant 
  $K^q_n(\LCO_q^{\eps_1},\LCO_q^{\eps_2},\dots,\LCO_q^{\eps_n})$ vanishes unless
  $\eps_1=\eps_2=\dots=\eps_{n-1}=*$ and $\eps_n=n-1$.
  In the latter case
  $$
  K^q_n(\LCO_q^*,\LCO_q^*,\dots,\LCO_q^*,\LCO_q^{n-1})
  =
  \rho(\LCO_q^* \LCO_q^* \dots \LCO_q^* \LCO_q^{n-1})
  = [n-1]_q!
  $$
\end{Theorem}

\begin{proof}
  This follows from \eqref{eq:GeneralCumulantExamples:qFock:annihilation}, namely
  $$
  \LCO^* e^{\ox k}
  = (1+q+\dots+q^{n-1})\, e^{\ox (k-1)}
  = [n]_q\, e^{\ox (k-1)}
  $$
  from which we infer that
  $$
  {\LCO^*}^{n-1} \LCO^{n-1} \Omega
  = [n-1]_q!\, \Omega
  .
  $$
  Now apply Proposition~\ref{prop:GeneralCumulantFock:Lukasiewiczcumulants}.
\end{proof}


In the case $q=0$ (free Fock space) we have a complete description of the
cumulants.
\begin{Proposition}
  \label{prop:GeneralCumulants:FreeFockCumulants}
  Let~$q=0$. Then for~$k_j\in\IZ$ the cumulants are
  $$
  K^\exchF_n(\LCO^{k_1}, \LCO^{k_2}, \dots,\LCO^{k_n})
  = \begin{cases}
      \rho(\LCO^{k_1} \LCO^{k_2} \cdots \LCO^{k_n}) = 1 & \text{if there is no non-trivial lattice subword} \\
      0 &\text{otherwise}
  \end{cases}
  .
  $$
  Here~$\LCO^{k}$ stands for~$(\LCO^*)^{-k}$ if~$k<0$.
\end{Proposition}
\begin{proof}
  If there is a lattice subword $\LCO^{k_p},\LCO^{k_{p+1}},\dots, \LCO^{k_q}$,
  then 
  $$
  \LCO^{k_p} \LCO^{k_{p+1}} \cdots \LCO^{k_q}
  = \rho(\LCO^{k_p} \LCO^{k_{p+1}} \cdots \LCO^{k_q}) I
  $$
  and
  $$
  (\LCO^{k_p})^\omega (\LCO^{k_{p+1}})^\omega \cdots (\LCO^{k_q})^\omega
  =\rho( (\LCO^{k_p})^\omega (\LCO^{k_{p+1}})^\omega \cdots (\LCO^{k_q})^\omega )
  = 0
  $$
If there is no such subword, we use the noncrossing moment-cumulant formula
$$
K^q_n(\LCO^{k_1},\LCO^{k_2},\dots, \LCO^{k_1})
= \sum_{\pi\in \NC_n}
   \rho(\LCO_{\pi(1)}^{k_1} \LCO_{\pi(2)}^{k_2}\cdots \LCO_{\pi(n)}^{k_n})
   \,
   \mu_{\NC}(\pi,\hat{1}_n)
$$
and the only partition which gives a nonzero contribution is $\pi=\hat{1}_n$.
\end{proof}
Note that for general $q$ this is no longer true.

Let us now compute the partitioned cumulants
of \L{}ukasiewicz words, i.e.\ 
words of the form $\LCO_{\pi(1)}^{\eps_1} \cdots \LCO_{\pi(n)}^{\eps_n}$
with $\eps_j\in\{*,0,1,2,\dots\}$.
By Proposition~\ref{prop:GeneralCumulantFock:PartitionedLukasiewiczcumulants}
the cumulants coincide with the expectations
$$
\rho(\LCO_{\pi(1)}^{\eps_1} \cdots \LCO_{\pi(n)}^{\eps_n})
$$
and we only need to compute those, in which each block $B=\{i_1 < i_2 < \dots < i_b\}$
satisfies~$\eps_{i_1}=\eps_{i_2}=\dots =\eps_{i_{b-1}}=*$
and $\eps_{i_b}=b-1$.
For noncrossing~$\pi$, it follows by pyramidal independence that
$$
\rho(\LCO_{\pi(1)}^{\eps_1} \cdots \LCO_{\pi(n)}^{\eps_n})
= \prod_B \rho(\prod_{i\in B} \LCO^{\eps_i})
;
$$
for more general $\pi$ we must count the right reduced crossings.
\begin{Corollary}
  \label{cor:qCumulantsN:partitionedcumulants}
  For a partition $\pi=\{B_1,B_2,\dots,B_p\}$ the corresponding
  cumulant of the $q$-Toeplitz operator is
  $$
  K^q_\pi(\LCO_q^* + \sum_{k=0}^\infty \frac{\alpha_{k+1}}{[k]_q!} \, \LCO_q^k)
  = q^{rrc(\pi)}\, \alpha_{\abs{B_1}} \alpha_{\abs{B_2}} \cdots \alpha_{\abs{B_p}}
  $$
  where $rrc(\pi)$ is the number of right reduced crossings
  $$
  rrc(\pi) = \#\{ (i<i'<j<j') : i,j\in B, i',j'\in B', j=\max B, j'=\max B' \}
  $$
\end{Corollary}

\begin{proof}
  Let us define an un-crossing map $\Phi:\Pi_n\to\Pi_n$.
  For $\pi\in\Pi_n$, sort its blocks according to their last elements.
  Let $B_0$ be the last block which is not an interval.
  Choose $j_2\in B_0$ maximal s.t.\ $j_2-1\not\in B_0$ and
  let $j_1\in B_0$ be its predecessor in $B_0$.
  In other words, $j_1$ and $j_2$ enclose the last ``hole'' in $B_0$.
  Let $\Phi(\pi)$ be the partition obtained by cyclically rotating the interval
  $(j_1,j_1+1,\dots,j_2-1)$ to $(j_1+1,j_1+2,\dots,j_2-1,j_1)$.
  Then $rrc(\pi) = rrc(\Phi(\pi)) + c$
  where $c$ is the number of right reduced arcs of $\pi$ which are crossed by
  the arc $(j_1,\max B_0)$.
  This number is equal to
  $$
  c = \sum_{\substack{B\in\pi\\ j_1<\max B < j_2}}
       \abs{B} - \abs{B\cap[j_1+1,j_2-1]}
  $$
  and it is also immediate that 
  \begin{multline*}
  \rho(\LCO_{\pi(1)}^{\eps_1}
       \LCO_{\pi(2)}^{\eps_2}
       \cdots
       \LCO_{\pi(j_1)}^*
       \LCO_{\pi(j_1+1)}^{\eps_{j_1+1}}
       \cdots
       \LCO_{\pi(j_2-1)}^{\eps_{j_2-1}}
       \cdots
       \LCO_{\pi(n)}^{\eps_{n}}
      )
\\
  = q^c\,
  \rho(\LCO_{\pi(1)}^{\eps_1}
       \LCO_{\pi(2)}^{\eps_2}
       \cdots
       \LCO_{\pi(j_1-1)}^{\eps_{j_1-1}}
       \LCO_{\pi(j_1+1)}^{\eps_{j_1+1}}
       \cdots
       \LCO_{\pi(j_2-1)}^{\eps_{j_2-1}}
       \LCO_{\pi(j_1)}^*
       \LCO_{\pi(j_2)}^{\eps_{j_2}}
       \cdots
       \LCO_{\pi(n)}^{\eps_{n}}
      )
  \end{multline*}
\end{proof}


\subsection{``Reduced'' $q$-cumulants}

\ddate{20.07.2001}
In this section we consider the cumulants found in
\cite{Anshelevich:2001:partitiondependent}%
\ppaper{Anshelevich}{Partition-dependent stochastic measures and $q$-deformed cumulants}{Preprint, 2001}
(see also \cite{Saitoh+Yoshida:2001:canonical})
which are weighted by another statistic on partitions,
the number of so-called reduced crossings.
Instead of taking powers of creation operators, one adds 
gauge operators to the scenery.

\begin{Definition}
  Let $T$ be an operator with dense domain $\alg{D}$.
  The \emph{gauge operator} $\gamma(T)$ on $\Fock_q(\hilb{H})$ with dense domain
  $\Fock_{alg}(\alg{D})$ is defined by
  $$
  \gamma(T)\,\Omega = 0
  \qquad
  \gamma(T)\, \eta_1 \oxfdots \eta_n
  = \sum q^{k-1} T\eta_k\oxf\eta_1 \oxfdots \hat{\eta}_k \oxfdots \eta_n
  $$
\end{Definition}
\begin{Proposition}[{\cite[Prop.~2.2]{Anshelevich:2001:partitiondependent}\paper[Prop.~2.2]{A}{}{}}]
  If $T$ is essential selfadjoint with dense domain $\alg{D}$ and $T(\alg{D})\subseteq\alg{D}$,
  then $\gamma(T)$ is essential selfadjoint with dense domain $\Fock_{alg}(\alg{D})$.
\end{Proposition}
Then consider processes of the form
$$
\gamma_I(\xi,T,\lambda) = \LCO_I(\xi) + \LCO_I^*(\xi) + \gamma_I(T) + \abs{I}\lambda
$$
where $\hilb{H}=L^2(\IR_+)\ox V$ is the underlying Hilbert space,
$\xi\in V$ is an analytic vector for $T:\alg{D}\subseteq V\to V$ and
$\LCO_I(\xi) = \LCO(\chi_I\ox\xi)$ etc.
In \cite[Prop.~2.2]{Anshelevich:2001:partitiondependent} stochastic measures
are used to find the cumulants of such processes.

This kind of independence can be reduced to the following symmetry.
Let $\hilb{H}=\ell_2\ox V$ and denote $\{e_j\}_{j\geq0}$ the canonical basis of $\ell_2$
and $P_j$ the one-dimensional projections on $e_j$.
For $v\in V$, $T\in B(V)$ (or densely defined), $\lambda\in\IC$ (or $\IR$ to be s.a.),
define
$$
\gamma(v,T,\lambda) = \LCO(e_0\ox v) + \LCO^*(e_0\ox v) + \gamma(P_0\ox T) + \lambda I
$$
Interchangeable copies are obtained by replacing $e_0$ with $e_j$:
$$
\gamma_j(v,T,\lambda) = \LCO_j(v) + \LCO_j^*(v) + \gamma_j(T) + \lambda
$$
where
$$
\LCO_j(v) = \LCO(e_j\ox v)
\qquad
\gamma_j(T) = \gamma(P_j\ox T)
$$

\begin{Remark}
  Note that all these operators are infinite divisible.
  This can be seen by embedding $l_2$ into $L_2(\IR)$
  by sending the basis elements $e_i$ to the characteristic
  function of unit intervals. Therefore the class of
  obtainable distributions is the same as in
  \cite {Anshelevich:2001:partitiondependent}.
  We chose to work with $l_2$ for the sake of simplicity.
\end{Remark}

We will consider the unital $*$-algebras $\alg{A}_i$ generated by 
$\alg{X}_i = \{\LCO_i(v), \LCO_i^*(v), \gamma_i(T) : v\in V, T\in B(V)\}$
and more generally, for an index set $I$, the algebra $A_I$ generated by $\alg{A}_i$, $i\in I$.
Then these algebras satisfy pyramidal independence.
First we need a lemma.
\begin{Lemma}[{\cite[Lemma~3.2]{Anshelevich:2001:partitiondependent}\paper[Lemma~3.2]{A}{}{}}]
  \label{lemma:qCumulants:Ansh:3.2}
  Let $E_j=\lin\{e_i\}_{i\in I_j}\subseteq \ell_2$, $j=1,2$ with $I_1\cap I_2=\emptyset$
  (i.e.\ $E_1\perp E_2$).
  Then we have
  \begin{enumerate}
   \item $\Fock_q(E_1\ox V)\ominus \IC\Omega \perp \Fock_q(E_2\ox V)\ominus \IC\Omega$
   \item Let $\eta\in\Fock_q(E_2\ox V)$, $X\in \alg{A}_1$,
    then $X\eta = (X-\rho(\eta)) \,\Omega\ox\eta + \rho(X)\,\eta$
  \end{enumerate}
\end{Lemma}
\begin{proof}
  The first part is clear since $E_1\ox V\perp E_2\ox V$.
  For the second part, observe that by orthogonality we have for $i\in I_1$
  \begin{align*}
    \LCO_i^*(v)\,\eta    &= 0 \\
    \LCO_i(v)\,\eta  &= (e_i\ox v) \ox \eta \\
    \gamma_i(T)\,\eta &= 0
  \end{align*}
  More generally, for $\eta_1\in\Fock_q(E_1\ox V)\ominus\IC_\Omega$,
  \begin{align*}
    \LCO_i^*(v)\,\eta_1\ox\eta   &= (\LCO_i^*(v)\,\eta_1)\ox\eta  \\
    \gamma_i(T)\,\eta_1\ox\eta &= (\gamma_i(T)\,\eta_1)\ox\eta 
  \end{align*}
  Thus $\eta$ is either unchanged or sent to $0$ and the claim follows.
\end{proof}
Pyramidal independence still holds.
\begin{Proposition}[{\cite[Lemma~3.3]{Anshelevich:2001:partitiondependent}\paper[Lemma~3.3]{A}{}{}}]
  Let $X,X'\in\alg{A}_{I}$, $Y\in \alg{A}_J$ with $I\cap J=\emptyset$,
  then $\rho(XYX') = \rho(XX')\,\rho(Y)$.
\end{Proposition}
\begin{proof}
  By the preceding lemma,
  \begin{align*}
    \langle \Omega, XYX'\Omega\rangle
    &= \langle X^*\Omega, \overset{\circ}{Y}\ox X'\Omega + \rho(Y)\, X'\Omega\rangle \\
    &= \rho(Y)\, \langle \Omega,XX'\Omega\rangle\\
    &= \rho(Y)\,\rho(XX')
  \end{align*}
\end{proof}

Now we want to compute cumulants of the generators $\LCO_i(v)$, $\LCO_i^*(v)$, and $\gamma_i(T)$.
In the following let $X_i$ denote one of $\LCO_j(v)$, $\LCO_j^*(v)$, $\gamma_j(T)$, $j\in I$,
$v\in V$, $T\in B(V)$.
First observe that the expectation of a word $X_1X_2\cdots X_n$ vanishes unless $X_1$ is an
annihilator and $X_n$ is a creator.
Using lemma~\ref{lemma:qCumulants:Ansh:3.2} again one sees the following more general fact.
Let $w=X_1X_2\cdots X_n$ be a word in generators, $X_j\in \alg{A}_{i_j}$
and denote $\pi\in\Pi_n$ the partition induced by the indices $i_j$.
Then again $\rho(X_1X_2\cdots X_n)=0$ unless each block of $\pi$ starts with an
annihilator and ends with a creator. Moreover the number of creators must equal the number
of annihilators in each block. The following proposition shows that the cumulants are
even more restrictive and behave like those of ``gaussian'' variables: only one pair of
creator/annihilator is allowed in each block.
\begin{Proposition}
  The cumulant $K^q_\pi(X_1,X_2,\dots,X_n)$ vanishes
  unless each block of $\pi$ starts with an annihilator, ends with a creator and otherwise
  only contains gauge operators.
  If these conditions are satisfied, the cumulant equals the expectation of the corresponding
  word.
\end{Proposition}
\begin{proof}
  We only give the proof for $\pi=\hat{1}_n$;
  it remains essentially the same in the general case.
  The cumulant is equal to the expectation of the ``discrete Fourier
  transform'' \externaleqref{I}{eq:GeneralCumulants:GoodFormula1}.
  If we abbreviate $e_\omega = \sum \omega^k e_k$ and
  $p_\omega = \sum \omega^k p_k$, we have
  \begin{align*}
    \LCO_0(v)^\omega &= \sum \omega^k \LCO(e_k\ox v) = \LCO(e_\omega\ox v) \\
    \LCO_0^*(v)^\omega   &= \LCO^*(e_{\bar{\omega}}\ox v) \\
    \gamma_0(T)^\omega &= \gamma(p_\omega\ox T)
  \end{align*}
  Now $p_\omega e_{\omega^m} = e_{\omega^{m+1}}$
  and 
  $$
  \langle e_{\bar{\omega}}, e_{\omega^m}\rangle
    = \sum \omega^k \omega^{mk} 
    = \sum \omega^{(m+1)k} = 0
  $$
  unless $m+1$ is a multiple of $n$.
  The action of $X_1^\omega X_2^\omega\cdots X_n^\omega$ on $\Omega$
  starts with a creator
  $$
  X_1^\omega X_2^\omega\cdots X_n^\omega\Omega
  = X_1^\omega X_2^\omega\cdots X_{n-1}^\omega (e_\omega\ox v)
  $$
  Then a mixture of creation, annihilation and gauge operators changes the 
  first component by either tensoring with $e_\omega$ ore multiplying with $p_\omega$.
  If there is an annihilator besides $X_1$, it encounters $e_{\omega^m}$ with $m<n-1$
  and all the inner products vanish.
  If there is a creator, it must be matched by an annihilator different from $X_1$,
  and again the inner products vanish.
  Hence, $X_1,X_2,\dots,X_n$ must have the structure claimed in the theorem.
  If this is the case, we have $X_1=\LCO_0^*(v_1)$, $X_n=\LCO_0(v_2)$,
  and $X_j=\gamma_0(T_j)$ for $2\leq j\leq n-1$, and
  \begin{align*}
    \rho(X_1^\omega X_2^\omega\cdots X_n^\omega)
    &= \langle e_{\bar{\omega}}\ox v_1, e_{\omega^{n-1}}\ox T_2 T_3 \cdots T_{n-1} v_2\rangle \\
    &= \langle e_{\bar{\omega}}, e_{\omega^{n-1}}\rangle
       \,
       \langle v_1, T_2 T_3 \cdots T_{n-1} v_2\rangle \\
    &= n\, \rho(X_1 X_2 \cdots X_n)
  \end{align*}
\end{proof}

\ddate{23.07.2001}

\begin{Definition}[{An un-crossing map \cite{Anshelevich:2001:partitiondependent}}]
On the set of partitions $\Pi_n$ define a map $\Phi:\Pi_n\to\Pi_n$
which fixes interval partitions and otherwise acts as follows.
Let $B$ be the last block which is not an interval
(if we sort blocks with respect to their maximal element).
Let $j_2=\max\{s\in B : s-1\not\in B\}$ (the start of the last subinterval of $B$)
and $j_1$ its predecessor in $B$.
These two numbers enclose a hole of $B$ and the map $\Phi$ moves this hole to the end of $B$:
Let $\alpha=((j_1+1)(j_1+2)\cdots b(B))^{b(B)-j_2+1}\in\SG_n$
and $\Phi(\pi) = \alpha(\pi)$, i.e.,
$i\sim_{\Phi(\pi)} j$
$\iff$
$\alpha^{-1}(i)\sim_\pi \alpha^{-1}(j)$.
We will need the number of blocks which end between $j_1$ and $j_2$ but do not start there,
$$
c_b(\pi) = \abs{\{ s : j_1 < b(B_s) < j_2  \}}
           -
           \abs{\{ s : j_1 < a(B_s) < j_2  \}}
$$
Then $rc(\pi) = rc(\Phi(\pi)) + c_b(\pi)$ and since iterating the map $\Phi$ ends at an interval
partition after at most $n$ steps, we have
$$
rc(\pi) = \sum_{k=0}^n c_b(\Phi^k(\pi))
.
$$
\end{Definition}

\begin{Theorem}[{\cite[Lemma~3.8]{Anshelevich:2001:partitiondependent}\paper[Lemma~3.8]{A}{}{}}]
  Let $X_j\in \alg{X}_{i_j}$ be generators and $\pi\in\Pi_n$ be the kernel
  of the index map~$j\mapsto i_j$, i.e., $p\sim_\pi q$ $\iff$ $i_p=i_q$.
  Assume that each block consists of gauge operators enclosed by an annihilator at
  the beginning and a creator at the end.
  Then
  $$
  \rho(X_1 X_2 \cdots X_n) = q^{rc(\pi)} \prod_{B\in\pi} \rho(X_B)
  $$
  that is, for such words $\rho$ is multiplicative modulo a factor $q^{rc(\pi)}$.
\end{Theorem}

\begin{proof}
  Let $B$ be the block containing $n$.
  If $B$ is an interval, we can factor it out by pyramidal independence.

  If not, consider the last block and let $j_1$, $j_2$ be as in the construction of the
  un-crossing map above.
  If this hole is a union of intervals, we can factor it out by pyramidal independence.

  Otherwise there are crossings.
  Let 
  $$
  \eta = X_{j_2}X_{j_2+1}\cdots X_n\Omega\in \hilb{H}^{\ox 1}
  $$
  and
  $$
  \xi = X_{j_1+1}X_{j_1+2}\cdots X_{j_2-1}\Omega\in \hilb{H}^{\ox c_b(\pi)}
  $$
  note that $X_{j_1}$ is either an annihilator or a gauge operator,
  therefore
  \begin{align*}
    X_{j_1}X_{j_1+1}\cdots X_n\Omega
    &= X_{j_1}X_{j_1+1}\cdots X_{j_2-1} \eta \\
    &= X_{j_1} (\xi\ox\eta)\\
    &= q^{c_b(\pi)} X_{j_1}\eta\ox \xi \\
    &= q^{c_b(\pi)}
       X_{j_1}X_{j_2}X_{j_2+1}\cdots X_n
       X_{j_1+1}X_{j_1+2}\cdots X_{j_2-1}
       \Omega
  \end{align*}
  because $c_b(\pi)$ is the number of creators between $X_{j_1}$ and $X_{j_2}$.
  Thus
  $$
  \rho(X_1X_2\cdots X_n)
  = q^{c_b(\pi)}
    \rho(X_1X_2\cdots X_{j_1}
         X_{j_2}X_{j_2+1}\cdots X_n
         X_{j_1+1}X_{j_1+2}\cdots X_{j_2-1}
        )
  $$
  and the partition determined by the permuted indices is exactly $\Phi(\pi)$.
\end{proof}

By multilinear expansion we get from this

\begin{Corollary}
  For $X_j=\LCO_0^*(v_j)+\LCO_0(v_j) + \gamma_0(T_j)+\lambda_j$
  the partitioned cumulants are multiplicative modulo a factor $q^{rc(\pi)}$:
  $$
  K^q_\pi(X_1 X_2 \cdots X_n) = q^{rc(\pi)} \prod_{B\in\pi} K^q_{\abs{B}}(X_B)
  $$
\end{Corollary}


%% file: SymmetricGroup.tex
\ssection{Fock spaces associated to characters of the infinite symmetric group}

Recently the calculations on another concrete Fock space have been carried out in
\cite{BozejkoGuta:2002:functors}. It was shown that using a certain embedding of pair partitions
into symmetric groups, one can evaluate characters of the latter and obtain
Fock states in this way.

\subsection{A simple case}

\begin{Definition}
  Let~$\pi\in\Pi^{(2)}_{2n}$ be a pair partition.
  There exists a unique noncrossing pair partition~$\hat\pi\in\NC^{(2)}_{2n}$
  such that the set of left points of the pairs in~$\pi$ and~$\hat\pi$ coincide.
  A \emph{cycle} in~$\pi$ is a sequence~$((l_1,r_1),\dots,(l_m,r_m))$ of pairs of~$\pi$
  such that the pairs $(l_1,r_2),(l_2,r_3),\dots,(l_m,r_1)$ belong to $\hat\pi$.
  The length of this cycle is~$m$. We denote by~$c(\pi)$ the number of cycles of~$\pi$
  and $c_m(\pi)$ the number of cycles of length~$m$.
\end{Definition}

Let $\pi=\{(l_1,r_1),\dots,(l_n,r_n)\}$ and
$\hat{\pi}=\{(l_1,\hat{r}_1),\dots,(l_n,\hat{r}_n)\}$
then define a permutation $\sigma$ by its images $\sigma(i)=j$ if $r_i=\hat{r}_j$,
i.e., $\hat{r}_i=r_{\sigma^{-1}(i)}$.
Then the number of cycles of $\pi$ is the number of cycles of $\sigma$.

Here the cumulant function which we want to consider is given by
$$
\mathbf{t}_N(\pi)  =\left(\frac{1}{N}\right)^{\abs{\pi}-c(\pi)}
$$
where $N\in\IZ\setminus\{0\}$ is a fixed integer.
The corresponding combinatorial Fock space can be realized as follows.
The function $\mathbf{t}_N$ actually comes from a character of the
infinite symmetric group which is given by
$$
\phi_N(\sigma) = \prod_{m\geq2} \left(\frac{1}{N}\right)^{(m-1)\,c_m(\sigma)}
$$

With this function there is associated the symmetrizator
\begin{align*}
  P_N^{(n)} 
  &=\sum_{\sigma\in\SG_n}
    \phi_N(\sigma)\, \tilde{U}^{(n)}(\sigma)\\
  &= \left(e+ \frac{1}{N} \tilde{U}^{(n)}((1,2))\right)
     \left(e + \frac{1}{N} (\tilde{U}^{(n)}((1,3))+\tilde{U}^{(n)}((2,3)))\right)
     \cdots \\
  &\phantom{=======} \hfil
     \cdots
     \left(e + \frac{1}{N}(
               \tilde{U}^{(n)}((1,n))
               +
               \tilde{U}^{(n)}((2,n))
               +
               \dots
               +\tilde{U}^{(n)}((n-1,n)))
     \right)
\end{align*}

On full Fock space $\Fock(\hilb{H})$ consider the deformed inner product
$$
\langle \xi, \eta\rangle_N = \langle \xi, P_N \eta \rangle
$$
i.e., for simple tensors the inner product is
$$
\langle \xi_1\oxfdots \xi_m, 
        \eta_1\oxfdots \eta_n\rangle_N
 = \delta_{mn}
   \sum_{\sigma\in\SG_n}
    \phi_N(\sigma)\, 
    \langle \xi_1,\eta_{\sigma^{-1}(1)} \rangle
    \langle \xi_2,\eta_{\sigma^{-1}(2)} \rangle
    \cdots
    \langle \xi_n,\eta_{\sigma^{-1}(n)} \rangle
$$
One can show that this is a positive bilinear form and after dividing through
its kernel and completing it we get a Hilbert space $\Fock_N(\hilb{H})$.
On these we have the following creation and annihilation operators
$$
  \LCO_N(\xi)\, \Omega = \xi
\qquad\qquad
  \LCO_N^*(\xi)\, \Omega = 0
$$
\begin{align*}
  \LCO_N(\xi)\, \xi_1\oxfdots \xi_n &= \xi \oxf \xi_1\oxfdots \xi_n \\ 
  \LCO_N^*(\xi)\, \xi_1\oxfdots \xi_n
  &= \langle \xi,\xi_1\rangle\,  \xi_2\oxfdots \xi_n \\ 
  &\phantom=
     +
     \frac{1}{N}
     \sum_{k=2}^n
     \langle \xi,\xi_k\rangle
     \xi_2\oxfdots\xi_{k-1} \oxf\xi_1 \oxf \xi_{k+1}\oxfdots \xi_n
\end{align*}
It follows that for a unit vector $\xi\in\hilb{H}$ we have the deformed factorial
function
$$
b_{n+1}
= \rho( L_N^*(\xi)^n L_N(\xi)^n )
= \left(
    1+\frac{1}{N}
  \right)
  \left(
    1+\frac{2}{N}
  \right)
  \cdots
  \left(
    1+\frac{n-1}{N}
  \right)
$$
and by Corollary~\ref{cor:GeneralCumulantFock:generalizedToeplitzcumulant} we
have the following formula for the cumulants.
\begin{Proposition}
  $$
  K^\exchm_n(L_N^* 
      +
      \sum_{k=1}^\infty 
       \frac{\alpha_k}%
            {\left(
               1+\frac{1}{N}
             \right)
             \left(
               1+\frac{2}{N}
             \right)
             \cdots
             \left(
               1+\frac{k-2}{N}
             \right)
           }
       L_N^{k-1})
  = \alpha_n
  $$
\end{Proposition}

\subsection{Some graph theory: the cycle cover polynomial}
\label{sec:SymmetricGroup:cyclecoverpolynomial}

In order to understand the partitioned moments and cumulants we need some graph theory.

Given a multidimensional Dyck word
$W=L_{i_1}^{\eps_1}L_{i_2}^{\eps_2}\dotsm L_{i_{2n}}^{\eps_{2n}}$ with $\eps_j\in\{*,1\}$,
the expectation
$$
\rho(W) = \sum_{\substack{\pi\in\Pi_{2n}^{(2)}\\ \pi\leq\underline{i}}} t^{n-c(\pi)}
$$
roughly has the following interpretation.
To the word $W$ we associate a digraph $\Gamma_W$.
First ignore the ``colors'' $i_j$ and
consider the onedimensional Dyck word $L_1^{\eps_1}L_1^{\eps_2}\dotsm L_1^{\eps_{2n}}$.
The number of $*$'s among the $\eps_j$ is the same as the number of $1$'s.
Let $\hat{\pi}\in \NC_{2n}^{(2)}$ be the unique noncrossing pair partition whose
left points are the indices $j$ with $\eps_j=*$.
We interpret this partition as a bijection $b$ between the annihilators and creators.
The graph $\Gamma_W$ has vertex set $V=\{1,\dots,n\}$ and we put an arrow
from $v_p$ to $v_{p'}$ if $L_{i_p}^*$ and $L_{i_{b(p')}}$ have the same color
and $b(p')>p$. 

In other words, we first construct a bipartite graph~$B$ on the two sets of vertices:
The set of annihilators~$V_*=\{p(1)<p(2)<\dots<p(n)\}$ 
and creators $V_1=\{q(1)<q(2)<\dots<q(n)\}$.
The edges are $E=\{\{p(r),q(s)\}: i_{p(r)}=i_{q(s)}\ \text{and}\ q(s)>p(r) \}$.
There is a unique noncrossing pair partition $\pi\in\NC_{2n}^{(2)}$
consisting of pairs $\{p<q\}$ s.t.\ $\eps_p=*$ and $\eps_q=1$.
It can be constructed recursively by connecting the rightmost annihilator
to its right neighbour, removing both from the word and repeating the procedure
on the new word. In particular, considering $\pi$ as a permutation (product of
transpositions), we have $\pi(p(n))=p(n)+1$.
Then
$$
\rho(L_{i_1}^{\eps_1}L_{i_2}^{\eps_2}\dotsm L_{i_{2n}}^{\eps_{2n}})
=\sum_{\mu} \frac{1}{N^{n-c(\sigma)}}
$$
where the sum runs over all complete matchings of the bipartite graph $B$
and $\sigma\in\SG_n$ is the permutation on $V_*$ which maps $p(i)$ to $\pi^{-1}\circ\mu(i)$,
if we consider $\pi$ and $\mu$ as bijections from~$V_*$ to~$V_1$.

Further we construct a directed graph $\Gamma$ with vertex set  $\{p(1),\dots,p(n)\}$
which we relabel as $v_1,\dots,v_n$ and we put a directed edge from $v_r$ to $v_s$
if $\{p(r),\pi(p(s))\}\in E(B)$, where $B$ is the bipartite graph constructed above.
This is the contraction of bipartite graph $B$ along the bijection $\pi$,
which was also constructed in \cite{Lass:2002:variations}.

In other words, $v_rv_s\in E(\Gamma)$ if and only if
$i_{p(r)}=i_{\pi(p(s))}=i_{q(\sigma(s))}$ and $\pi(p(s))>p(r)$.
A complete matching of~$B$ corresponds to a partition of~$\Gamma$ into
hamiltonian cycles, that is, a partition of the vertices into vertex-disjoint cycles,
and the number of cycles of~$\sigma$ is equal to the number of components
of the hamiltonian partition (also called a cycle cover or $2$-factor).

We are therefore led to the so called \emph{cycle cover polynomial} $C_c(\Gamma;x)$.
This is a specialization of both the cover polynomial $C!(\Gamma;x,y)$ of Chung and Graham 
\cite{ChungGraham:1995:cover} and the geometric cover polynomial $C(\Gamma;x,y)$ of
D'Antona and Munarini~\cite{DAntonaMunarini:2000:cyclepath}, 
which have been proposed as a kind of Tutte polynomial for digraphs.
Namely we have 
\begin{align*}
  C!(\Gamma;x,y) &= \sum_{(C,P)} x^{\abs{C}} y^{\underline{\abs{P}}}\\
  C(\Gamma;x,y) &= \sum_{(C,P)} x^{\abs{C}} y^{\abs{P}}\\
  C_c(\Gamma;x) &= C!(\Gamma;x,0) = C(\Gamma;x,0) = \sum_C x^{\abs{C}}
\end{align*}
here the sum runs over all cycle-path covers $(C,P)$ (cycle covers $C$, respectively) 
of $\Gamma$ with $\abs{C}$ cycles and $\abs{P}$ paths.
Loops are interpreted as cycles of length $1$ and single vertices as paths
of length $0$ and
$y^{\underline{n}}=y(y-1)\dotsm(y-n+1)$ is the falling factorial function.

\begin{Theorem}
$$
\rho(L_{i_1}^{\eps_1}L_{i_2}^{\eps_2}\dotsm L_{i_{2n}}^{\eps_{2n}})
= \frac{1}{N^n}\, C_c(\Gamma;N)
$$
\end{Theorem}

This can be seen explicitly as discussed above or by observing that both
quantities satisfy the same recursive relation:

The creation operators satisfy the following ``commutation relation'':
\begin{equation}
  \label{eq:SymmetricGroup:simplecommutationrelation}
  L_i^*L_j = \delta(i,j)I + \frac{1}{N}\,\dG_{ji}
\end{equation}
with $\dG_{ji}=\dG(\ket{e_j}\,\bra{e_i})$
where for an operator~$T\in B(\hilb{H})$ we denote by~$\dG(T)$ 
its differential second quantization
$$
\dG(T)\, \xi_1\oxfdots\xi_n
= \sum_{i=1}^n
   \xi_1\oxfdots T\xi_i\oxfdots\xi_n
.
$$
These operators act like derivations on analytic vectors:
$$
\dG_{ji} L_{i_1}\dotsm L_{i_n}\Omega
= \sum_{k=1}^n \delta(i,i_k)\, L_{i_1}\dotsm L_{i_{k-1}} L_j L_{i_{k+1}}\dotsm L_{i_n}\Omega
$$
i.e., it replaces~$L_i$ by~$L_j$ and for any noncommutative polynomial~$P$
\begin{equation}
  \label{eq:SymmetricGroup:simpledGamma}
  \dG_{ji}P(L_1,L_2,\dots)\,\Omega
  = D_i P(L_1,L_2,\dots)[L_j]\,\Omega,
\end{equation}
the noncommutative derivative of $P$ with respect to $L_i$ in the direction $L_j$..
In particular $\dG_{ii} P\Omega = n_i P\Omega$, where $n_i$ is the total degree
of $L_i$ in $P$.

Therefore we have the following recursive procedure to compute the expectation 
of a multidimensional Dyck word
$W=L_{i_1}^{\eps_1}L_{i_2}^{\eps_2}\dotsm L_{i_{2n}}^{\eps_{2n}}$,
namely to look for the rightmost annihilator $L_{i_{p(n)}}^*$, replace
$L_{i_{p(n)}}^* L_{i_{p(n)+1}}$ by 
$\delta(i_{p(n)},i_{p(n)+1})I + \frac{1}{N}\,\dG_{i_{p(n)+1},i_{p(n)}}$
and apply~\eqref{eq:SymmetricGroup:simpledGamma}.

In terms of cover polynomials this corresponds to the cut and fuse recursion,
which is a convenient way to compute the cover polynomial.
Fix an edge $e=v_1v_2$ of $\Gamma$. Then it is easy to see that
\begin{equation}
  \label{eq:SymmetricGroup:coverpolynomialrecursion}
  C(\Gamma;x,y)
  = \begin{cases}
      x\,C(\Gamma\setminus v;x,y) + C(\Gamma\setminus e;x,y) &\text{if $e$ is a loop on the vertex $v_1=v_2=v$}\\
      C(\Gamma\setminus e;x,y) + C(\Gamma/e;x,y) &\text{if $e$ is not a loop}
    \end{cases}
\end{equation}

where the quotient $\Gamma/e$ is obtained from $\Gamma$ by fusing $v_1=e_+$
and $v_2=e_-$ along $e$, i.e., in the new graph the vertices $v_1$ and $v_2$ are identified
and all edges with source $v_1$ or target $v_2$ are deleted.
As initial condition we put $C(\Gamma_n;x,y)=y^n$ for the graph $\Gamma_n$ on $n$ vertices
with no edges and $C(\Gamma_0;x,y)=1$ for the empty graph $\Gamma_0$.
The polynomial of Chung and Graham satisfies exactly the same recursion, but with initial 
condition $C!(\Gamma_n;x,y)=y^{\underline{n}}$.
We are only interested in the cycle cover polynomial, for which the powers of $y$ are not relevant.

\begin{Claim}
  There is a one-to-one correspondence between the recursion
  $$
  \rho(L_{i_1}^*\dotsm L_{i_k}^*L_{i_{k+1}}\dotsm L_{i_{2n}})
  = \rho(L_{i_1}^*\dotsm L_{i_{k-1}}^{\eps_{k-1}}
  (\delta(i_k,i_{k+1})\,I+\frac{1}{N}\,\dG_{i_{k+1},i_k})L_{i_{k+2}}\dotsm L_{i_{2n}}) 
  $$
  and~\eqref{eq:SymmetricGroup:coverpolynomialrecursion}.
\end{Claim}
There are two cases to consider, depending on whether $i_{p(n)}=i_{p(n)+1}$ or not.
First assume~$i_{p(n)}=i_{p(n)+1}=:i$, then with $k=p(n)$ we have
\begin{align*}
  \rho(L_{i_1}^*\dotsm L_{i_k}^*L_{i_{k+1}}\dotsm L_{i_{2n}})
  &= \rho(L_{i_1}^*\dotsm L_{i_{k-1}}^{\eps_{k-1}} (I+\frac{1}{N}\,\dG_{ii})L_{i_{k+2}}\dotsm L_{i_{2n}}) \\
  &= (1+\frac{n_i}{N})
     \rho(L_{i_1}^*\dotsm L_{i_{k-1}}^{\eps_{k-1}} L_{i_{k+2}}\dotsm L_{i_{2n}}) \\
\end{align*}
where $n_i$ is the total degree of $L_i$ in $L_{i_{k+2}}\dotsm L_{i_{2n}}$.
On the graph side, there is a loop from $v_n$ to itself, and therefore removing
successively all edges emanating from~$v_n$ we obtain
\begin{align*}
  C_c(\Gamma;x)
  &= x\,C_c(\Gamma\setminus v_n;x) + C_c(\Gamma\setminus v_nv_n;x)\\
  &= x\,C_c(\Gamma\setminus v_n;x)
     +
     \sum_{\substack{t\\ i_p(n)=i_{\pi(p(t))}\\ \pi(p(t))>p(n)}}
     C_c(\Gamma/ v_nv_t;x)
\end{align*}
The second sum runs over all edges emanating from $v_n$.
The remaining term $C_c(\Gamma\setminus\{v_nv_n,v_nv_{t_1},\dots,v_nv_{t_k}\};x)$
vanishes because there is no cycle through $v_n$ anymore.

Now we show that each $\Gamma/v_nv_t$ is isomorphic to $\Gamma\setminus v_n$.
Indeed $v_rv_s\in E(\Gamma\setminus v_n)$ if and only if $v_rv_s\in E(\Gamma)$,
while $v_rv_s\in E(\Gamma/v_nv_t)$ if and only if either
$s\ne t$ and $v_rv_s\in E(\Gamma)$
or
$s=t$ and $v_rv_n\in E(\Gamma)$.
The first case being trivial, consider the case $s=t$.
Then $v_nv_t\in E(\Gamma)$ implies $i_{p(n)}=i_{\pi(p(t))}$ and
$\pi(p(t))>p(n)$.
Now $v_rv_t\in E(\Gamma)$ if and only if $i_{p(r)}=i_{\pi(p(t))}$ and $\pi(p(t))>p(r)$.
The second condition is redundant because already $p(n)>p(r)$.

On the other hand, 
$$
v_rv_n\in E(\Gamma) \iff i_{p(r)}=i_{\pi(p(n))}\ \text{and}\ \pi(p(n))=p(n)+1>p(r)
$$
and again the second condition is automatically satisfied.
Moreover $i_{\pi(p(n))}=i_{p(n)+1}=i_{p(n)}=i_{\pi(p(t))}$.
Thus $v_rv_t\in E(\Gamma)$ if and only if $v_rv_n\in E(\Gamma)$ and
finally
$v_rv_s\in E(\Gamma/v_nv_t)$ if and only if $v_rv_s\in E(\Gamma)$ for all $s$.

Now assume that $i_{p(n)}\ne i_{p(n)+1}$.
Then with $k=p(n)$
\begin{align*}
  \rho(L_{i_1}^*\dotsm L_{i_k}^*L_{i_{k+1}}\dotsm L_{i_{2n}})
  &= \rho(L_{i_1}^*\dotsm L_{i_{k-1}}^{\eps_{k-1}} \frac{1}{N}\dG_{i_{k+1},i_k}L_{i_{k+2}}\dotsm L_{i_{2n}}) \\
  &= \frac{1}{N}
     \sum_{r=k+2}^{2n}
      \delta(i_k,i_r)\,
      \rho(L_{i_1}^*\dotsm L_{i_{k-1}}^{\eps_{k-1}} L_{i_{k+2}}\dotsm
           L_{i_{r-1}} L_{i_{k+1}} L_{i_{r+1}}\dotsm L_{i_{2n}}) \\
\end{align*}
On the graph side, this corresponds to the identity
$$
C_c(\Gamma;x)
= \sum_{e_-=v_n} C(\Gamma/e;x)
= \sum_{\substack{t\\ i_{p(n)}=i_{\pi(p(t))}\\ \pi(p(t))>p(n)}}
   C_c(\Gamma/v_nv_t;x)
$$
and there is a one to one correspondence between the summands:
Let $r=\pi(p(t))>p(n)$,
then the graph $\Gamma_r$ corresponding to 
$\rho(L_{i_1}^*\dotsm L_{i_{k-1}}^{\eps_{k-1}} L_{i_{k+2}}\dotsm
L_{i_{r-1}} L_{i_{k+1}} L_{i_{r+1}}\dotsm L_{i_{2n}})$
has vertices $v_1,\dots,v_{n-1}$, partition $\pi'=\pi\setminus\{p(n),p(n)+1\}$
and edges 
\begin{align*}
  v_xv_y\in E(\Gamma_r)
  &\iff \begin{cases}
          \text{$i_{p(x)}=i_{\pi(p(y))}$ and $\pi'(p(y))>p(x)$} & \text{if $y\ne t$}\\
          \text{$i_{p(x)}=i_{p(n)+1}$} & \text{if $y=t$}
        \end{cases}\\
  &\iff \begin{cases}
          v_xv_y\in E(\Gamma)& \text{if $y\ne t$}\\
          v_xv_n \in E(\Gamma) & \text{if $y=t$}
        \end{cases}\\
  &\iff v_xv_y\in E(\Gamma/v_nv_t)
\end{align*}

\begin{Remark}
  Not every digraph is the digraph of a Dyck word.
  A necessary condition is, that the vertices with common successors can be
  linearly ordered. Denote by $R(v)$ the set of successors of a vertex $v$.
  Given two vertices $v_i$ and $v_j$, if there is a vertex $v_k$
  such that both $R(v_i)\cap R(v_k)\ne\emptyset$ and $R(v_j)\cap
  R(v_k)\ne\emptyset$,
  then either $R(v_i)\subseteq R(v_j)$ or $R(v_j)\subseteq R(v_i)$.
\end{Remark}

\subsection{The Vershik-Kerov construction}

We recall the construction of the representations of Vershik and Kerov
\cite{VershikKerov:1981:characters} which give rise to the irreducible characters
of the infinite symmetric group $\SG_\infty$
which were first found by Thoma \cite{Thoma:1964:unzerlegbaren},
namely
$$
\oldphi_{\alpha,\beta}(\sigma)
= \prod_{m\geq2}
   \left(
     \sum_{i=1}^\infty \alpha_i^m + (-1)^{m+1} \sum_{i=1}^\infty \beta_i^m
   \right)^{\rho_m(\sigma)}
$$
where $\alpha_1\geq\alpha_2\geq\dots\geq0$ and $\beta_1\geq\beta_2\geq\dots\geq0$
are given sequences such that $\sum\alpha_i+\sum\beta_i\leq1$
and $\rho_m(\sigma)$ is the number of cycles of length~$m$ of~$\sigma$.

To keep notations simple we consider only the the case where $\beta_i=0$ and $\sum\alpha_i=1$.
In this case $\alpha_n$ can be interpreted as a probability measure~$m^{(\alpha)}$ on~$\IN$ with
$m^{(\alpha)}(\{k\})=\alpha_k$.
We consider the elements of the space $\XX_n=\IN^n$ with the product measure~$m_n^{(\alpha)}$
as words~$x=x_1x_2\dotsm x_n$ and abbreviate the value of the measure by
$$
\alpha_x= m_n^{(\alpha)}(\{x\})=\prod_{i=1}^n \alpha_{x_i}.
$$
$\SG_n$ acts on $\XX_n$ by $(\sigma x)_i= x_{\sigma^{-1}(i)}$. 
We write $x\sim y$ if $x$ and $y$ are in the same orbit, i.e.,
if there exists a permutation $\sigma\in\SG_n$ such that $y=\sigma x$.
Now consider the space of ``rearrangements''
$$
\tilde{\XX}_n = \{(x,y)\in \XX_n\times\XX_n : x\sim y\}
.
$$
If we equip the space of functions $f:\tilde{\XX}\to\IC$
with the Hilbert norm
$$
\norm{f}_2^2 = \int_{\XX_n} \sum_{y\sim x} \abs{f(x,y)}^2\,dm_n^{(\alpha)}(x) 
= \sum_{x\in\XX_n} \alpha_x \sum_{y\sim x} \abs{f(x,y)}^2
,
$$
then the Hilbert space
$$
V_n^{(\alpha)} 
= \{f:\tilde{\XX}_n\to\IC : \norm{f}_2^2 < \infty \}
$$
carries a unitary representation $U^{(\alpha)}_n$ of $\SG_n$
$$
(U_n^{(\alpha)}(\sigma)\,h)(x,y) = h(\sigma^{-1}x,y)
$$
This space has an orthogonal basis $\{\bb_{x,y}: (x,y)\in\tilde{\XX}_n\}$
of functions $\bb_{x,y}(x',y')=\delta(x,x')\,\delta(y,y')$.
The action of $\SG_n$ is
$$
U_n^{(\alpha)}(\sigma)\, \bb_{x,y} = \bb_{\sigma x,y}
$$
and their inner products are
$$
\langle \bb_{x',y'},\bb_{x'',y''} \rangle_{V_n^{(\alpha)}}
= \delta(x',x'')\,\delta(y',y'')\,\alpha_{x'}
$$

\begin{Theorem}[\cite{VershikKerov:1981:characters}]
  Let $1_n=\sum_{x\in\tilde{\XX}_n} \bb_{x,x}$ be the diagonal,
  then
  $$
  \langle U_n^{(\alpha)}(\sigma)\, 1_n,1_n \rangle
  = m_n^{(\alpha)}(\{x : \sigma x=x \})
  = \oldphi_{\alpha,0}(\sigma)
  $$
\end{Theorem}
\begin{Proposition}
  The maps $j_n:V_n^{(\alpha)}\to V_{n+1}^{(\alpha)}$ defined by
  $$
  j_n h(x,y) = \delta(x_{n+1},y_{n+1})\,h(x_1\dotsm x_n,y_1\dotsm y_n),
  $$
  i.e.,
  $$
  j_n \bb_{x,y} = \sum_{z\in\XX_1} \bb_{xz,yz}
  $$
  clearly 
  satisfy the intertwining relation~\eqref{eq:GeneralCumulantFock:intertwining}
  and their adjoints are given by
  $$
  j_n^* \bb_{xz,yz'} = \delta(z,z')\, \alpha_z \bb_{x,y}
  .
  $$
\end{Proposition}

We can therefore construct the symmetric Fock space spanned by the vectors
$$
\bb_{x,y} \ox_s \xi_1\oxf\xi_2\oxfdots\xi_n
= \frac{1}{n!}
  \sum_{\sigma\in\SG_n}
   \bb_{\sigma x,y}\ox\xi_{\sigma^{-1}(1)}
                   \oxf\xi_{\sigma^{-1}(2)}
                   \oxfdots\xi_{\sigma^{-1}(n)}
$$
It will be notationally more convenient to work with the right creation
and annihilation operators \eqref{eq:GeneralCumulantFock:creationoperatordefinition}
\begin{align*}
  R(\xi_{n+1})\, \bb_{x,y} \ox_s \xi_1\oxf\xi_2\oxfdots\xi_n
  &= (n+1) \,j_n\bb_{x,y} \ox_s \xi_1\oxf\xi_2\oxfdots\xi_{n+1}\\
  &= (n+1)\sum_{z\in\XX_1} \bb_{xz,yz} \ox_s \xi_1\oxf\xi_2\oxfdots\xi_{n+1}\\
\end{align*}
The annihilation operator is given by 
\eqref{eq:GeneralCumulantFock:annihilationoperator} and can be simplified to
\begin{align*}
  R^*(\xi)\,&\bb_{x,y} \ox_s \xi_1\oxf\xi_2\oxfdots\xi_{n+1}\\
  &= \frac{1}{n+1}
     \sum_{i=1}^{n+1}
      \langle \xi_i,\xi \rangle \,
      j_n^*U_{n+1}(\tau_{i,n+1})\, \bb_{x,y}\ox_s \xi_1\oxfdots\xi_{i-1}\oxf\xi_{n+1}\oxf\xi_{i+1}\oxfdots\xi_n\\
  &= \frac{1}{n+1}
     \sum_{i=1}^{n+1}\,
      \langle \xi_i,\xi \rangle \,
      \delta(x_i,y_{n+1})\,
      \alpha_{x_i}\,
      \bb_{x_1\dotsm \Check{x}_i\dotsm x_n,y_1\dotsm y_n}\ox_s \xi_1\oxfdots\Check{\xi}_i\oxfdots\xi_n\\
\end{align*}
\begin{Definition}
  For an operator $T\in B(\hilb{H})$ we
  define weighted differential second quantisation operators
  $$
  \dG^{(k)}(T)\,\bb_{x,y}\ox_s\xi_1\oxfdots\xi_n
  = \sum_{i=1}^n \alpha_{x_i}^k\, \bb_{x,y}\ox_s\xi_1\oxfdots\xi_{i-1}\oxf T\xi_i\oxf\xi_{i+1}\oxfdots\xi_n
  $$
\end{Definition}
Then we have the following ``commutation relation''.
\begin{Lemma}
  \begin{equation}
    \label{eq:SymmetricGroup:commutationrelation}
    R^*(\xi) R(\eta)
    = \langle \eta,\xi \rangle\,I+\dG^{(1)}(\ket{\eta}\,\bra{\xi})
  \end{equation}
  In particular for standard basis vectors $\xi=e_i$ and $\eta=e_j$ we have
  \begin{equation}
    \label{eq:SymmetricGroup:commutationrelationij}
    R_i^*R_j = \delta(i,j)\,I+\dG_{ji}^{(1)}
  \end{equation}
  where $\dG_{ji}^{(k)} = \dG^{(k)}(\ket{e_j}\,\bra{e_i})$.
\end{Lemma}
\begin{proof}
  \begin{align*}
    R^*(\xi) & R(\xi_{n+1})\,\bb_{x,y}\ox_s\xi_1\oxfdots\xi_n\\
    &= (n+1)\,
       R^*(\xi)
       \sum_{z\in\XX_1} 
        \bb_{xz,yz}
        \ox_s\xi_1
        \oxfdots\xi_{n+1}\\
    &= \sum_{z\in\XX_1} 
        \sum_{i=1}^{n+1}
         \langle \xi_i,\xi\rangle\,
         j_n^*U_{n+1}(\tau_{i,n+1})\,\bb_{xz,yz}
         \ox_s\tilde{U}_{n+1}(\tau_{i,n+1})\,\xi_1
         \oxfdots\xi_{n+1}\\
    &= \sum_{z\in\XX_1} 
        \sum_{i=1}^{n}
         \langle \xi_i,\xi\rangle\,
         \delta(x_i,z)\,
         \alpha_z\,
         \bb_{x_1\dotsm \Check{x}_i\dotsm x_nz,y}
         \ox_s\xi_1
         \oxfdots\Check{\xi}_i
         \oxfdots\xi_{n+1}\\
    &\phantom{==}   +
       \sum_{z\in\XX_1} 
        \langle \xi_{n+1},\xi\rangle\,
         \alpha_z\,
         \bb_{x,y}
         \ox_s\xi_1
         \oxfdots\xi_n\\
    &= \sum_{i=1}^{n}
        \langle \xi_i,\xi\rangle\,
        \delta(x_i,z)\,
        \alpha_{x_i}\,
        \bb_{x_1\dotsm \Check{x}_i\dotsm x_nx_i,y}
        \ox_s\xi_1
        \oxfdots\Check{\xi}_i
        \oxfdots\xi_{n+1}\\
    &\phantom{==}   +
        \langle \xi_{n+1},\xi\rangle\,
         \bb_{x,y}
         \ox_s\xi_1
         \oxfdots\xi_n\\
    &= \sum_{i=1}^{n}
        \langle \xi_i,\xi\rangle\,
        \delta(x_i,z)\,
        \alpha_{x_i}\,
        \bb_{x,y}
        \ox_s\xi_1
        \oxfdots\xi_{i-1}
        \oxf\xi_{n+1}
        \oxf\xi_{i+1}
        \oxfdots\xi_n\\
    &\phantom{==}   +
        \langle \xi_{n+1},\xi\rangle\,
         \bb_{x,y}
         \ox_s\xi_1
         \oxfdots\xi_n
  \end{align*}
\end{proof}

It will also be convenient to have weighted creation and annihilation operators
at our disposal.
\begin{Definition}
  For $k\in\IN_0$ let
  $$
  R^{(k)}(\xi)\,\bb_{x,y}\ox_s \xi_1\oxfdots\xi_n
  = (n+1)\sum_{z\in\XX_1} \alpha_z^k\,\bb_{xz,yz}\ox_s\xi_1\oxfdots\xi_n\oxf\xi
  $$
  in particular, $R^{(0)}(\xi)=R(\xi)$.
\end{Definition}
Then the action of $\dG_{ji}$ on analytic vectors is like a derivation:
\begin{Lemma}
  $$
  \dG_{ji}^{(k)} R_{i_n}^{(k_n)}\dotsm R_{i_1}^{(k_1)}\,\Omega 
  = \sum_{j=1}^n 
     \delta(i,i_j)\,
     R_{i_n}^{(k_n)}\dotsm R_{i_{j+1}}^{(k_{j+1})} R_j^{(k_j+k)}  R_{i_{j-1}}^{(k_{j-1})} \dotsm R_{i_1}^{(k_1)}\,\Omega
  $$
\end{Lemma}
The adjoint of a weighted creation operator is a weighted annihilation operator,
as expected.
\begin{Proposition}
  $$
  R^{(k)}(\xi)^*\, \bb_{x,y}\ox_s \xi_1\oxfdots\xi_n
  = \frac{1}{n}
    \sum_{i=1}^n 
    \delta(x_i,y_n)\,
    \alpha_{x_i}^{k+1}\,
    \langle \xi_i,\xi \rangle\,
    \bb_{x_1\dotsm \Check{x}_i\dotsm x_n,y_1\dotsm y_{n-1}}
    \ox_s \xi_1\oxfdots\Check{\xi}_i\oxfdots\xi_n
  $$
\end{Proposition}

The general ``commutation relation'' now is as follows.
\begin{Proposition}
  $$
  R^{(k)*}(\eta)R^{(m)}(\xi)
  = \langle\xi,\eta\rangle
    \sum_z \alpha_z^{m+k+1} I
    +\dG^{(m+k+1)}(\ket{\xi}\,\bra{\eta})
  $$
  and in particular
  $$
  R_i^{(k)*} R_j^{(m)}
  = \delta(i,j)\sum_z \alpha_z^{m+k+1} I 
    +\dG_{ji}^{(m+k+1)}
  $$
\end{Proposition}

In this case the expectations of multidimensional Dyck words is given
by evaluations of the \emph{cycle indicator polynomial}
of D'Antona and Munarini.

\subsection{More graph theory: the cycle indicator polynomial}

\begin{Definition}[{\cite{DAntonaMunarini:2000:cyclepath}}]
  Let $\Gamma$ be a digraph.
  The \emph{cycle-path indicator polynomial} of $\Gamma$ is the multivariate
  polynomial
  $$
  I(\Gamma;\bx,\by)
  = \sum_{(C,P)}
     \prod_{\gamma\in C} x_{\abs{\gamma}} 
     \prod_{\pi\in P} y_{\abs{\pi}} 
  $$
  where the sum runs over all cycle-path coverings~$(C,P)$ of $\Gamma$,
  and $\abs{\gamma}$ (resp.\ $\abs{\pi}$) denotes the number of vertices
  of a cycle $\gamma$ (resp.\ a path $\pi$).
\end{Definition}
The cycle-path indicator polynomial satisfies a recursion similar to the
recursion \eqref{eq:SymmetricGroup:coverpolynomialrecursion}
of the geometric cover polynomial.
However in order to do the cut and fuse operations, one somehow has to remember
the deleted vertices. This is done by putting weights on the vertices.
\begin{Definition}
  Let~$w:V(\gamma)\to\IN$ be a weight function on the vertices of the
  digraph~$\Gamma$ and denote the weighted graph by~$\Gamma_w=(\Gamma,w)$.
  For a path $\pi=v_1\dotsm v_k$ let $w(\pi)=w(v_1)+\dots+w(v_k)$
  and similarly the weight of a cycle is the sum of the weights of its
  vertices.
  Define the cycle-path indicator polynomial of $\Gamma_w$ as
  $$
  I(\Gamma_w;\bx,\by)
  = \sum_{(C,P)}
     \prod_{\gamma\in C} x_{w(\gamma)}
     \prod_{\pi\in P} y_{w(\pi)}
  $$
  An unweighted graph therefore has the same cycle-path indicator polynomial
  as the same graph with weight~$1$ on each vertex.

  We define the cut and fuse operations on weighted digraphs as follows.
  \begin{align*}
    \Gamma_w\setminus v &= (\Gamma\setminus v,w|_{\Gamma\setminus v}) &\text{for $v\in V(\Gamma)$}\\
    \Gamma_w\setminus e &= (\Gamma\setminus e,w)                      &\text{for $w\in E(\Gamma)$}\\
    \Gamma_w/e &= (\Gamma/e,w^*)                      &\text{for $e=v_1v_2\in E(\Gamma)$}\\
  \end{align*}
  where we denote $v_0$ the collapsed vertex and
  $$
  w^*(v)=
  \begin{cases}
    w(v) & v\ne v_0\\
    w(v_1)+w(v_2) & v=v_0
  \end{cases}
  $$
\end{Definition}

\begin{Proposition}[{\cite{DAntonaMunarini:2000:cyclepath}}]
  The cycle-path indicator polynomial satisfies the following cut-and-fuse
  recursion. Fix an edge~$e$, then
  $$
  I(\Gamma_w;\bx,\by) 
  = \begin{cases}
      x_{w(v)} I(\Gamma_w\setminus v;\bx,\by) + I(\Gamma_w\setminus e;\bx,\by) 
      &\text{if $e=vv$ is a loop}\\
      I(\Gamma_w\setminus e;\bx,\by) + I(\Gamma_w/ e;\bx,\by) 
      &\text{if $e=v_1v_2$ is not a loop}
  \end{cases}
  $$
  If $\Gamma$ has no edges, then
  $$
  I(\Gamma_w;\bx,by) = \prod_v y_{w(v)}
  $$
  and if $\Gamma$ is empty, then $I(\emptyset;\bx,\by)=1$.
\end{Proposition}
We are only interested in the \emph{cycle indicator polynomial}
$$
I_c(\Gamma_w;\bx)=I(\Gamma_w;\bx,\boldsymbol{0}).
$$

Again the commutation relation~\eqref{eq:SymmetricGroup:commutationrelationij}
corresponds to cut and fuse operations and we obtain the following formula.
Roughly the expectation of a weighted multidimensional Dyck word is given 
by the cycle indicator of the weighted graph $\Gamma_w$, where $\Gamma$
is constructed as in section~\ref{sec:SymmetricGroup:cyclecoverpolynomial} and
the weight of a vertex is computed from the weight of the corresponding
annihilator plus one plus
the weight of the creator to which it is connected by $\pi$.
\begin{Theorem}
  $$
  \rho(R_{i_1}^{(k_1)\eps_1}\dotsm R_{i_1}^{(k_{2n})\eps_{2n}})
  = I_c(\Gamma_w;\bx)
  $$
  where $\Gamma_w$ is the graph with weights 
  $w(v_i)=k_{p(i)}+1 + k_{\pi(p(i))}$
  and $x_1=1$,
  $$
  x_k = \sum_{i=1}^\infty \alpha_i^k + (-1)^{k+1} \sum_{i=1}^\infty \beta_i^k
  $$
  for $k\geq2$.
\end{Theorem}

\begin{Theorem}
  Any positive definite functional on pair partitions whose values only depends on
  the graphs~$\Gamma(\pi)$ is given by a (not necessarily irreducible) character
  on the infinite symmetric group.
  The multiplicative ones are exactly those which come from irreducible characters
  of the symmetric group.
\end{Theorem}
\begin{proof}
  It suffices to show this for pair partitions.
  We will construct the representation of~$\SG_n$ as in the proof 
  of~\cite[Theorem~2.7]{GutaMaassen:2002:generalised}.
  Let~$H_n$ be the space spanned by products
  $$
  R_{\sigma(1)}  R_{\sigma(2)}\dotsm R_{\sigma(n)}\Omega
  $$
  this space carries a natural representation~$U_n$ of $\SG_n$, namely
  $$
  U_n(\sigma)  R_{\tau(1)}  R_{\tau(2)}\dotsm R_{\tau(n)}\Omega
  =  R_{\sigma\tau(1)}  R_{\sigma\tau(2)}\dotsm R_{\sigma\tau(n)}\Omega
  $$
  and therefore for any partition with given cycle structure the expectation function 
  is given by
  $$
  \langle U_n(\sigma) R_n R_{n-1}\dotsm R_1\Omega, R_n R_{n-1}\dotsm R_1\Omega \rangle
  =  \langle
       R_1^* R_2^*\dotsm R_n^*  R_{\sigma(n)} R_{\sigma(n-1)}\dotsm R_\sigma(1)\Omega,
       \Omega
     \rangle
  = t(\pi)
  $$
  where $\sigma$ has the same cycle structure as $\pi$.
\end{proof}


%% file: GeneralCumulant3.bbl
\providecommand{\bysame}{\leavevmode\hbox to3em{\hrulefill}\thinspace}
\begin{thebibliography}{vLM96}

\bibitem[Ans01]{Anshelevich:2001:partitiondependent}
Anshelevich, M., \emph{Partition-dependent stochastic measures and $q$-deformed
  cumulants}, Doc. Math. \textbf{6} (2001), 343--384.

\bibitem[BG02]{BozejkoGuta:2002:functors}
Bo{\.z}ejko, M., and Guta, M., \emph{Functors of white noise associated to
  characters of the infinite symmetric group}, Comm. Math. Phys. \textbf{229}
  (2002), 209--227.

\bibitem[BKS97]{BozejkoKuemmererSpeicher:1997:qGaussian}
Bo{\.z}ejko, M., K{\"u}mmerer, B., and Speicher, R., \emph{$q$-{G}aussian
  processes: non-commutative and classical aspects}, Comm. Math. Phys.
  \textbf{185} (1997), 129--154.

\bibitem[BS94]{BozejkoSpeicher:1994:completely}
Bo{\.z}ejko, M., and Speicher, R., \emph{Completely positive maps on {C}oxeter
  groups, deformed commutation relations, and operator spaces}, Math. Ann.
  \textbf{300} (1994), 97--120.

\bibitem[CG95]{ChungGraham:1995:cover}
Chung, F. R.~K., and Graham, R.~L., \emph{On the cover polynomial of a
  digraph}, J. Combin. Theory Ser. B \textbf{65} (1995), 273--290.

\bibitem[DM00]{DAntonaMunarini:2000:cyclepath}
D'Antona, O.~M., and Munarini, E., \emph{The cycle-path indicator polynomial of
  a digraph}, Adv. in Appl. Math. \textbf{25} (2000), 41--56.

\bibitem[FB70]{FrischBourret:1970:parastochastics}
Frisch, U., and Bourret, R., \emph{Parastochastics}, J. Mathematical Phys.
  \textbf{11} (1970), 364--390.

\bibitem[GM02]{GutaMaassen:2002:generalised}
Gu{\c{t}}{\u{a}}, M., and Maassen, H., \emph{Generalised {B}rownian motion and
  second quantisation}, J. Funct. Anal. \textbf{191} (2002), 241--275.

\bibitem[Gut01]{Guta:2001:qproduct}
Guta, M., \emph{The $q$-product of generalised {B}rownian motions}, Preprint,
  September 2001.

\bibitem[Las02]{Lass:2002:variations}
Lass, B., \emph{Variations sur le th\`eme {$E+\overline E=XY$}}, Adv. in Appl.
  Math. \textbf{29} (2002), 215--242.

\bibitem[Leh03]{Lehner:2002:Cumulants2}
Lehner, F., \emph{Cumulants in noncommutative probability theory {II}.
  {G}eneralized {G}aussian random variables}, Probab. Theory Related Fields
  \textbf{127} (2003), 407--422, arXiv:math.CO/0210443.

\bibitem[Leh04]{Lehner:2002:Cumulants1}
Lehner, F., \emph{Cumulants in noncommutative probability theory {I}.
  {N}oncommutative exchangeability systems}, Math. Zeitschr. \textbf{to appear}
  (2004), arXiv:math.CO/0210442.

\bibitem[MN97]{MingoNica:1997:crossings}
Mingo, J.~A., and Nica, A., \emph{Crossings of set-partitions and addition of
  graded-independent random variables}, Internat. J. Math. \textbf{8} (1997),
  645--664.

\bibitem[Nic95]{Nica:1995:oneparameter}
Nica, A., \emph{A one-parameter family of transforms, linearizing convolution
  laws for probability distributions}, Comm. Math. Phys. \textbf{168} (1995),
  187--207.

\bibitem[Nic96]{Nica_1996:crossings}
Nica, A., \emph{Crossings and embracings of set-partitions and $q$-analogues of
  the logarithm of the {F}ourier transform}, Discrete Math. \textbf{157}
  (1996), 285--309.

\bibitem[SY01]{Saitoh+Yoshida:2001:canonical}
Saitoh, N., and Yoshida, H., \emph{A canonical random variable for the
  $q$-deformed moments-cumulants formula}, Preprint, 2001.

\bibitem[Tho64]{Thoma:1964:unzerlegbaren}
Thoma, E., \emph{Die unzerlegbaren, positiv-definiten {K}lassenfunktionen der
  abz{\"a}hlbar unendlichen, symmetrischen {G}ruppe}, Math. Z. \textbf{85}
  (1964), 40--61.

\bibitem[VK81]{VershikKerov:1981:characters}
Vershik, A.~M., and Kerov, S.~V., \emph{Characters and factor representations
  of the infinite symmetric group}, Dokl. Akad. Nauk SSSR \textbf{257} (1981),
  no.~5, 1037--1040.

\bibitem[vLM96]{vanLeeuwenMaassen:1996:obstruction}
van Leeuwen, H., and Maassen, H., \emph{An obstruction for $q$-deformation of
  the convolution product}, J. Phys. A \textbf{29} (1996), 4741--4748.

\bibitem[Voi85]{Voiculescu:1985:symmetries}
Voiculescu, D., \emph{Symmetries of some reduced free product ${C}\sp
  \ast$-algebras}, Operator algebras and their connections with topology and
  ergodic theory (Bu\c steni, 1983), Springer, Berlin, 1985, pp.~556--588.

\bibitem[Voi86]{Voiculescu:1986:addition}
Voiculescu, D., \emph{Addition of certain noncommuting random variables}, J.
  Funct. Anal. \textbf{66} (1986), 323--346.

\end{thebibliography}
